\documentclass[12pt]{amsart}
\usepackage{amsmath}
\usepackage{amsfonts}
\usepackage{amssymb}
\usepackage[all]{xy}           %xypic macro for latex2.09
\usepackage{bm}
\usepackage{bbding}
\usepackage{txfonts}
\usepackage{amscd}
\usepackage{xspace}
\usepackage[shortlabels]{enumitem}
\newtheorem{theorem}{Theorem}[section]

\newtheorem{lemma}{Lemma}[section]
\newtheorem{prop}{Proposition}[section]
\newtheorem{cor}{Corollary}[section]
\newtheorem{definition}{Definition}[section]
\newtheorem{example}{Example}[section]

\numberwithin{equation}{section}
\newcommand{\Hom}{\mathrm{Hom}}

\newcommand{\g}{\mathfrak g}
\newcommand{\gl}{\frak{gl}}
\newcommand{\der}{{\rm Der}}
\newcommand{\F}{\mathbb{F}}

\newcommand{\dM}{\mathrm{d}}
\newcommand{\frkX}{\mathfrak X}
\newcommand{\huaV}{\mathcal{V}}
\newcommand{\LD}{\mathrm{LieDer}}
\newcommand {\emptycomment}[1]{}
\begin{document}

\title[]{Cohomologies of $3$-Lie algebras with derivations}

\author{Senrong Xu}
\address{Xu: School of Mathematical Sciences, Jiangsu University, Zhenjiang 212013, Jiangsu, China}
\email{senrongxu@ujs.edu.cn}

\author{Jiefeng Liu$^{\ast}$}
\address{Liu: School of Mathematics and Statistics, Northeast Normal University, Changchun 130024, Jilin, China}
\email{liujf534@nenu.edu.cn}
\thanks{$^{\ast}$ the corresponding author}

\vspace{-5mm}

%\date{\today}

\begin{abstract}
In this paper, we consider a $3$-Lie algebra with a derivation (called a {\rm $3$-LieDer} pair). We define cohomology for a {\rm $3$-LieDer} pair with coefficients in a representation.  We use this cohomology to study deformations and abelian extensions of {\rm $3$-LieDer} pairs. We give the notion of a {\rm $3$-Lie$2$Der} pair, which can be viewed as the categorification of a  {\rm $3$-LieDer} pair. We show that skeletal {\rm $3$-Lie$2$Der} pairs  are classified by triples given by $3$-LieDer pairs, representations and $3$-cocycles.  We define crossed modules of {\rm $3$-LieDer} pairs and show that there exists a one-to-one correspondence between strict {\rm $3$-Lie$2$Der} pairs and crossed modules of {\rm $3$-LieDer} pairs.

\end{abstract}



\maketitle

\tableofcontents

\allowdisplaybreaks

%\end{document}

\section{Introduction}
Generalizations of Lie algebras to higher arities, including 3-Lie algebras and more generally, n-Lie algebras (also called Filippov algebras) (\cite{Fi}), have attracted attention from both mathematics and physics. $n$-Lie algebras are the algebraic structure corresponding to Nambu mechanics (\cite{Na}). In particular, the structure of $3$-Lie algebras is applied to the study of the supersymmetry and gauge symmetry transformations of the word-volume theory of multiple M2-branes (\cite{BaLa08,BaLa09}) and the Bagger-Lambert theory has a novel local gauge symmetry which is based on the metric $3$-Lie algebras (\cite{MeFi08,MeFi09}). Many extensive literatures are related to this pioneering work. See \cite{CS,HHM,P} and the review article \cite{deAI} for more details. Structure theory of $n$-Lie algebras are widely studied. Representation theory of $n$-Lie algebras was studied by Kasymov  in \cite{Ka} and cohomology theory of $n$-Lie algebras was studied by Takhtajan and Gautheron in \cite{Ga,Ta2}. Deformations of  $n$-Lie algebras were studied in \cite{ABM,F OARRILL,Ta2}. See \cite{deAI,LMS,Xu} for more details on extensions of $n$-Lie algebras.

Derivations are useful tool to study various algebraic structures. Derivations can be used to construct homotopy Lie algebras (\cite{Vo}), deformation formulas (\cite{CGG}) and differential
Galois theory (\cite{Mag}). They also play an essential role in  mathematical logic (\cite{RS}), functional analysis (\cite{BCS,BdS}) and control theory (\cite{AK,AT}). In \cite{DL,Lo}, the authors study algebras with derivations from the operadic point of view. Recently, the authors in \cite{TFS} study the cohomology, extensions and deformations of Lie algebras with derivations (called LieDer pairs). The results of \cite{TFS} have been extended to
associative algebras with derivations (called AssDer pairs) and  Leibniz algebras with derivations (called LeibDer pairs) in \cite{Da1,Da2}.

In this paper, we consider  $3$-Lie algebras with derivations (called {\rm $3$-LieDer} pairs). We develop a cohomology theory of {\rm $3$-LieDer} pairs that controls the extensions and deformations of {\rm $3$-LieDer} pairs.

The paper is organized as follows. In Section \ref{s2}, we introduce representations  and cohomologies of {\rm $3$-LieDer} pairs. We construct a degree $-1$ graded Lie algebra structure on the underlying graded vector space of the cochain complex such that the Maurer-Cartan elements of this graded Lie algebra are precisely {\rm $3$-LieDer} pairs. In Section \ref{s3}, we consider deformations of {\rm $3$-LieDer} pairs. We define one-parameter formal deformations of {\rm $3$-LieDer} pairs and show that $n$-infinitesimals of equivalent one-parameter formal deformations of a {\rm $3$-LieDer} pair are in the same cohomology class. We study deformations of order $n$ of {\rm $3$-LieDer} pairs and show that the obstruction of a deformation of order $n$ extending to a deformation of order $n+1$ can be controlled by the third cohomology groups. In Section \ref{abel}, we study abelian extensions of {\rm $3$-LieDer} pairs and show that abelian extensions of {\rm $3$-LieDer} pairs are classified by the second cohomology groups. In Section \ref{s5}, we introduce the concept of {\rm $3$-Lie$2$Der} pairs, which is the categorification of {\rm $3$-LieDer} pairs.  We show that there is a one-to-one correspondence between skeletal {\rm $3$-Lie$2$Der} pairs and the triples consisted of {\rm $3$-LieDer} pairs, representations and $3$-cocycles. We define crossed modules of {\rm $3$-LieDer} pairs and show that there exists a one-to-one correspondence between strict {\rm $3$-Lie$2$Der} pairs and crossed modules of {\rm $3$-LieDer} pairs. Furthermore, we show that there exists a one-to-one correspondence between strict isomorphism classes of strict {\rm $3$-Lie$2$Der} pairs and equivalent classes of crossed modules of {\rm $3$-LieDer} pairs.

In this paper, we work over an algebraically closed field $\F$ of characteristic $0$ and all the vector spaces are over $\F$ and finite-dimensional.

\noindent
%{\bf Acknowledgements. }  We give warmest thanks to---

\section{Representations and cohomologies of {\rm $3$-LieDer} pairs}\label{s2}
In this section, we consider  {\rm $3$-LieDer} pairs. We give  representations and cohomologies of   {\rm $3$-LieDer} pairs. Finally, we construct a degree $-1$ graded Lie algebra  whose Maurer-Cartan elements are precisely {\rm $3$-LieDer} pairs.

\subsection{$3$-$\LD$ pairs and their representations}
\begin{definition}
  A {\bf 3-Lie algebra}  is a vector space $L$ together with a skew-symmetric linear map  $[\cdot,\cdot,\cdot]_L:
\otimes^3 L \rightarrow L$ such that the following Fundamental Identity  holds:
\begin{eqnarray}
\label{fundamental}&&[x_1,x_2,[x_3,x_4,x_5]_L]_L\\
\nonumber&=&[[x_1,x_2,x_3]_L,x_4,x_5]_L+[x_3,[x_1,x_2,x_4]_L,x_5]_L+[x_3,x_4,[x_1,x_2,x_5]_L]_L,
\end{eqnarray}
where $x_1,x_2,x_3,x_4,x_5\in L$.
\end{definition}

Elements in $\wedge^2 L$ are called {\bf fundamental objects} of the $3$-Lie algebra $(L,[\cdot,\cdot,\cdot]_L)$. There is a bilinear operation $[\cdot,\cdot]_F$ on $\wedge^2 L$, which is given by
\begin{eqnarray}
[X,Y]_F=[x_1,x_2,y_1]_L\wedge y_2+y_1 \wedge[x_1,x_2,y_2]_L, \quad \forall X=x_1\wedge x_2, ~Y=y_1\wedge y_2\in\wedge^2 L.\label{courant}
\end{eqnarray}
It is well known that $(\wedge^2 L,[\cdot,\cdot]_F)$ is a Leibniz algebra, which plays an important role in the theory of $3$-Lie algebras. The Leibniz rule can be written as
\begin{eqnarray}
[X,[Y,Z]_F]_F=[[X,Y]_F,Z]_F+[Y,[X,Z]_F]_F, \quad \forall X, Y , Z\in \wedge^2 L.
\end{eqnarray}
Moreover, the fundamental identity (\ref{fundamental}) is equivalent to
\begin{eqnarray}
[X,[Y,z]]-[Y,[X,z]]=[[X,Y]_F,z], \quad \forall X, Y \in \wedge^2 L, ~z \in L.\label{fun}
\end{eqnarray}

\begin{definition}\label{3lierep}{\rm(\cite{Ka})}.
	A  {\bf representation} of a $3$-Lie algebra $L$ on a vector space $V$ is a  linear map  $\rho:\wedge^{2}L \rightarrow {\rm End}(V)$ such that for all $x_{1},x_{2},x_{3},x_{4}\in L$,
	\begin{align*}
	\rho(x_{1},x_{2})\rho(x_{3},x_{4})\!&=\!\rho([x_{1},x_{2},x_{3}]_L,\!x_{4})\!+\!\rho(x_{3},[x_{1},x_{2},x_{4}]_L)\!+\!\rho(x_{3},x_{4})\rho(x_{1},x_{2}),\\
	\rho(x_{1},[x_{2},x_{3},x_{4}]_L)\!&=\!\rho(x_{3},x_{4})\rho(x_{1},x_{2})\!-\!\rho(x_{2},x_{4})\rho(x_{1},x_{3})\!+\!\rho(x_{2},x_{3})\rho(x_{1},x_{4}).
	\end{align*}
\end{definition}
In the sequel, a representation is denoted by $(V;\rho).$
\begin{prop}\label{semidirect}
Let $L$ be a $3$-Lie algebra, $V$ a vector space and $\rho:\wedge^2 L \to \gl(V)$ a bilinear map. Then $(V;\rho)$ is a representation of $L$ on $V$ if and only if there is a $3$-Lie algebra structure on the direct sum $L \oplus V$ of vector spaces (called the {\bf semidirect product}) defined by
$$[x_1+v_1,x_2+v_2,x_3+v_3]_\rho=[x_1,x_2,x_3]_L+\rho(x_1,x_2)v_3+\rho(x_3,x_1)v_2+\rho(x_2,x_3)v_1,~ \forall x_i \in L, v_i \in V.$$
A semidirect product $3$-Lie algebra is denoted by $L \ltimes_\rho V$ or simply $L \ltimes V$.
\end{prop}

The  cohomology theory for a $3$-Lie algebra $L$ with coefficients in a representation $(V;\rho)$ is given as follows. Denote by $C^n(L;V):=\Hom(\otimes^{n-1}(\wedge^2 L)\wedge L,V)$, the space of $n$-cochains. The corresponding coboundary operator $\dM: C^n(L;V) \to C^{n+1}(L;V)$ is defined by
 \begin{eqnarray}
 ~\nonumber &&(\dM f_n)(\frkX_1,...,\frkX_n,z)\\
 ~\nonumber&=& \sum\limits_{1 \le j < k\le n} {{{( - 1)}^j}f_n({\frkX_1},\cdots,\mathop {\hat {{\frkX_j}}  ,\cdots,{\frkX_{k - 1}},[\frkX_j,\frkX_k]_F,\cdots,\frkX_n,z)} } \\
 ~ \nonumber&& + \sum\limits_{j=1}^n(-1)^j{f_n(\frkX_1,\cdots, {\hat {\frkX_j}},\cdots,\frkX_n,[x_j,y_j,z]_L)}\\
 ~ \nonumber&& +\sum\limits_{j=1}^n{(-1)^{j+1}\rho(\frkX_j)f_n(\frkX_1,\cdots,{\hat {{\frkX_j}}} ,\cdots,\frkX_n,z)}\\
 ~ && +(-1)^{n+1}\rho(y_n,z)f_n(\frkX_1,\cdots,\frkX_{n-1},x_n)+(-1)^n\rho(x_n,z)f_n(\frkX_1,\cdots,\frkX_{n-1},y_n)\label{coboundary}
 \end{eqnarray}
 for all $f_n \in C^n(L;V)$ and  $\frkX_i=x_i\wedge y_i \in \wedge^2 L,~z \in L,~i=1,2,\cdots,n$, where $[\cdot,\cdot]_F$ is given by \eqref{courant}. We denote the corresponding $n$-th cohomology group by $H^{n}(L;V)$.

Recall that a {\bf derivation} on a $3$-Lie algebra $L$ is a linear map $\theta_L:L\to L$ satisfying
\begin{equation}
  \theta_L[x,y,z]_L=[\theta_L(x),y,z]_L+[x,\theta_L(y),z]_L+[x,y,\theta_L(z)]_L,\quad\forall~x,y,z\in L.
\end{equation}
We denote the set of derivations on the $3$-Lie algebra $L$ by $\der(L)$.

\begin{definition}
A {\bf $3$-LieDer pair} consists of a $3$-Lie algebra $(L,[\cdot,\cdot,\cdot]_L)$ and a derivation $\theta_L\in \der(L)$. We denote a $3$-$\LD$ pair by $(L,[\cdot,\cdot,\cdot]_L,\theta_L)$ or simply by $(L,\theta_L)$.
\end{definition}

%Next we give an example of 3-Lie algebra derivation which is induced by Lie algebra and its derivation.
\begin{example}
Let $L$ be the $4$-dimensional $3$-Lie algebra with basis $\{e_1,e_2,e_3,e_4\}$ and bracket
$$[e_1,e_2,e_3]_L=e_4, [e_1,e_2,e_4]_L=e_3, [e_1,e_3,e_4]_L=e_2, [e_2,e_3,e_4]_L=e_1.$$
Then any derivation $D=(d_{ij})\in\der(L)$ has the form as follows
$$D\begin{pmatrix}
	e_1\\
	e_2\\
e_3\\
e_4\end{pmatrix}=\begin{pmatrix}
	0&d_{12}&d_{13}&d_{14}\\
	d_{12}&0&d_{23}&d_{24}\\
-d_{13}&d_{23}&0&d_{34}\\
d_{14}&-d_{24}&d_{34}&0
	\end{pmatrix}\begin{pmatrix}
	e_1\\
	e_2\\
e_3\\
e_4
	\end{pmatrix}.$$
\end{example}
See \cite{AEM} for more examples and applications of derivations on $3$-Lie algebras.

\emptycomment{\begin{example}(\cite{AEM}).
Let $\g$ be a Lie algebra and $D$ is a Lie algebra derivation on $\g$. Set $\tau\in\g^*$. Define a 3-ary operation on $\g$ by (denote it by $\g_\tau$)
$$[x,y,z]\triangleq\mathop{\circlearrowleft}\limits_{{\rm x}, {\rm y},{\rm z}}\tau({\rm x})[{\rm y},{\rm z}],$$
where $\tau\in\g^*$ and
$$\mathop{\circlearrowleft}\limits_{{\rm x}, {\rm y},{\rm z}}\tau({\rm x})[{\rm y},{\rm z}]=\tau(x)[y,z]+\tau(y)[z,x]+\tau(z)[x,y].$$
Suppose that $\tau$ is a trace function on $\g$. Namely, $\tau([\g,\g])=0$. Then $\tau\circ D$ is also a trace function and hence $\g_\tau$ and $\g_{\tau\circ D}$ are 3-Lie algebras. By a direct check we get $D$ is a 3-Lie algebra derivation of $\g_\tau$ if and only if $\g_{\tau\circ D}$ is abelian.
\end{example}}

\begin{definition}\label{3liepmor}
	Let $(L,[\cdot,\cdot,\cdot]_L,\theta_L)$ and $(K,[\cdot,\cdot,\cdot]_K,\theta_K)$ be two {\rm 3-LieDer} pairs. A {\rm 3-LieDer} pair {\bf morphism} from $(L,[\cdot,\cdot,\cdot]_L,\theta_L)$  to $(K,[\cdot,\cdot,\cdot]_K,\theta_K)$ is a $3$-Lie algebra morphism $\eta:L\rightarrow K$ such that $\eta\circ\theta_{L}=\theta_{K}\circ\eta$.
\end{definition}

\emptycomment{By \cite[Lemma 2.3]{LMS}, we see that $L\ltimes V$ is a 3-Lie algebra. Any derivation of $L\ltimes V$ has the form $\left(\begin{array}{cc}f_{11}&f_{12}\\f_{21}&f_{22}\\\end{array}\right)$, where $f_{11}\in {\rm End}(L)$, $f_{12}\in {\rm Hom}(V, L)$, $f_{21}\in {\rm Hom}(L, V)$ and $f_{22}\in {\rm End}(V)$. Derivations of the form $\left(\begin{array}{cc}f_{11}&0\\0&f_{22}\\\end{array}\right)$ will be denoted as $(f_{11}, f_{22})\in {\rm End}(L)\times {\rm End}(V)$ for brevity.}

Next we give the definition of representations of {\rm $3$-LieDer} pairs.
\begin{definition}\label{reppa}
A {\bf representation} of a {\rm $3$-LieDer pair} $(L,\theta_L)$ on a vector space $V$ is a pair $(\rho,\theta_V)$, where $\rho:\wedge^2L\rightarrow {\rm End}(V)$ and $\theta_V\in{\rm End}(V)$, such that $\rho$ is a representation of the $3$-Lie algebra $L$ on $V$ and satisfies
\begin{align}
	\theta_V\circ\rho(x,y)-\rho(x,y)\circ\theta_V=\rho(\theta_L(x),y)+\rho(x,\theta_L(y)).\label{prf}
	\end{align}
 We denote a representation of a {\rm $3$-LieDer pair} by $(V;\rho,\theta_V)$.
\end{definition}

\begin{example}
	Let $(L,[\cdot,\cdot,\cdot]_L,\theta_L)$ be a $3$-$\LD$ pair. Then $(L;ad,\theta_L)$ is a representation of the $3$-$\LD$ pair $(L;\theta_L)$, which is called the {\bf adjoint representation} of  $(L,\theta_L)$.
\end{example}

\begin{lemma}\label{3liederiff}
Let $L$ be a $3$-Lie algebra and $(V;\rho)$ a representation of $L$. Then $(V;\rho,\theta_V)$ is a representation of a $3$-$\LD$ pair  $(L,\theta_L)$ if and only if   $(L\ltimes_\rho{V},\theta_L+\theta_V)$ is a $3$-$\LD$ pair, which we call the {\bf semi-direct product} of the $3$-$\LD$ pair $(L,\theta_L)$ by the representation $(V;\rho,\theta_V)$.
\end{lemma}
{\bf Proof.}
Assume that  $(L\ltimes_\rho{V},\theta_L+\theta_V)$ is a $3$-$\LD$ pair. Since $L\ltimes_\rho{V}$ is a semi-direct product $3$-Lie algebra, $\rho$ is a representation of the $3$-Lie algebra $L$ on $V$. Since $\theta_L+\theta_V$ is a derivation on the semi-direct product $3$-Lie algebra $L\ltimes_\rho{V}$, for $x,y\in L, v\in V$, we have
\begin{eqnarray*}
0&=&  (\theta_L+\theta_V)([x,y,v]_{\rho})-([(\theta_L+\theta_V)(x),y,v]_{\rho}+[x,(\theta_L+\theta_V)(y),v]_{\rho}+[x,y,(\theta_L+\theta_V)(v)]_{\rho})\\
&=&(\theta_V\circ\rho(x,y))(v)-(\rho(\theta_L(x),y)+\rho(x,\theta_L(y))+\rho(x,y)\circ\theta_V)(v),
\end{eqnarray*}
which implies that
$$\theta_V\circ\rho(x,y)=\rho(\theta_L(x),y)+\rho(x,\theta_L(y))+\rho(x,y)\circ\theta_V.$$
Thus $(V;\rho,\theta_V)$ is a representation of a $3$-$\LD$ pair  $(L,\theta_L)$.

The converse can be proved similarly. We omit the details.
 \hfill $\Box$

\vspace{3mm}

Let $(V;\rho,\theta_V)$ be a representation of a $3$-$\LD$ pair  $(L,\theta_L)$. With a slight abuse of notation, we still use $\theta_L$ to denote the endomorphism on $\wedge^2L$ given by
\begin{eqnarray*}
	\theta_L(X)=\theta_L(x)\wedge y+x\wedge\theta_L(y),\ \forall X=x\wedge y\in\wedge^2L.
\end{eqnarray*}
By a direct calculation, we have
\begin{align}
\theta_L([X,Y]_F)=[\theta_L(X),Y]_F+[X,\theta_L(Y)]_F,\label{txy}
\end{align}
where $[\cdot,\cdot]_F$ is given by \eqref{courant} and $X,Y\in \wedge^2L$. Thus $\theta_L$ is a derivation on the Leibniz algebra $(L,[\cdot,\cdot]_F)$.

\subsection{Cohomologies of $3$-$\LD$ pairs}
Let $(V;\rho,\theta_V)$ be a representation of a $3$-$\LD$ pair  $(L,\theta_L)$. We define cochain groups by $\mathcal{C}_{\rm 3-LieDer}^1(L;V):={\rm Hom}(L,V)$ and
\begin{eqnarray}
\mathcal{C}^{p}_{{\rm 3-LieDer}}(L;V):=C^p(L;V)\times C^{p-1}(L;V),\quad p\geq2.
\end{eqnarray}

For $p\geq1$, we  define an operator $\delta: C^p(L;V)\rightarrow C^p(L;V)$ by
\begin{align}
&&&(\delta\alpha_p)(X_1\otimes \cdots\otimes X_{p-1}\otimes z)\notag\\
&&=&\sum\limits_{i=1}^{p-1}\alpha_p(X_1\otimes\cdots\otimes\theta_L(X_i)\otimes\cdots X_{p-1}\otimes z)\notag\\
&&&+\alpha_p(X_1\otimes \cdots\otimes X_{p-1},\theta_L(z))-\theta_V(\alpha_p(X_1\otimes \cdots\otimes X_{p-1}\otimes z)),\label{delta}
\end{align}
where $X_{i}\in \wedge^{2}L$ and $z\in L$.

Define $\partial: \mathcal{C}^{1}_{{\rm 3-LieDer}}(L;V)\rightarrow \mathcal{C}^{2}_{{\rm 3-LieDer}}(L;V)$ by
\begin{equation}
  \partial \alpha_1=(\dM \alpha_1,-\delta \alpha_1),\quad \alpha_1\in {\rm Hom}(L,V).
\end{equation}

For $p\geq 2$, define $\partial: \mathcal{C}^{p}_{{\rm 3-LieDer}}(L;V)\rightarrow \mathcal{C}^{p+1}_{{\rm 3-LieDer}}(L;V)$ by
\begin{eqnarray}
\partial(\alpha_p,\beta_{p-1})=({\rm d}\alpha_p,{\rm d}\beta_{p-1}+(-1)^p\delta\alpha_p).
\end{eqnarray}

\begin{lemma}\label{dd=dd}
The operator $\delta$ and $\dM$ are commutative, i.e.,${\rm d}\circ\delta=\delta\circ {\rm d}$.
\end{lemma}
{\bf Proof.} It follows by a straightforward tedious calculations. \hfill $\Box$

\vspace{3mm}

\begin{prop}\label{pp=0}
The map $\partial$ is a coboundary operator, i.e. $\partial\circ\partial=0$.
\end{prop}
{\bf Proof.} Let $(\alpha_p,\beta_{p-1})\in \mathcal{C}^{p}_{{\rm 3-LieDer}}(L;V)$. Then, by Lemma \ref{dd=dd}, we get
\begin{align*}
(\partial\circ\partial)(\alpha_p,\beta_{p-1})&=\partial({\rm d}\alpha_p,{\rm d}\beta_{p-1}+(-1)^p\delta\alpha_p)\\
&=({\rm d}({\rm d}\alpha_p),{\rm d}({\rm d}\beta_{p-1}+(-1)^p\delta\alpha_p)+(-1)^{p+1}\delta({\rm d}\alpha_p))\\
&=(0,0+(-1)^pd(\delta\alpha_p)-(-1)^p\delta({\rm d}\alpha_p))\\
&=(0,0).
\end{align*}\hfill $\Box$

\begin{definition}
 Let $(V;\rho,\theta_V)$ be a representation of a {\rm 3-LieDer} pair $(L,\theta_L)$.   The cohomology of the cochain complex $(\mathcal{C}^{\ast}_{{\rm 3-LieDer}}(L;V),\partial)$ is taken to be the {\bf cohomology of the $3$-$\LD$ pair $(L;\theta_L)$}.
 The corresponding $p$-th cohomology group is defined by
\begin{align}
\mathcal{H}_{\rm 3-LieDer}^p(L;V)=\mathcal{Z}_{\rm 3-LieDer}^p(L;V)/\mathcal{B}_{\rm 3-LieDer}^p(L;V),
\end{align}
where $\mathcal{Z}_{\rm 3-LieDer}^p(L;V)$ is the space of $p$-cocycles and $\mathcal{B}_{\rm 3-LieDer}^p(L;V)$ is the space of $p$-coboundaries.
 \end{definition}
%In the following, we use the symbol $\partial_{\rm ad}$ to refer  the coboundary operator associated to the
%adjoint representation of a {\rm 3-LieDer} pair.

In particular, when $p=1$, we have the following
\begin{cor}
	Let $(V;\rho,\theta_V)$ be a representation of a {\rm 3-LieDer} pair $(L,\theta_L)$. Then
	\begin{align*}
	\mathcal{H}_{\rm 3-LieDer}^1(L;V)=\{\alpha\in{\rm Hom}(L,V)\ |\ \alpha[x,y,z]_L&=\rho(x,y)\alpha(z)+\rho(y,z)\alpha(x)-\rho(x,z)\alpha(y),\\ \alpha\circ\theta_L&=\theta_V\circ\alpha,\quad \forall~x,y,z\in L \}.
	\end{align*}
\end{cor}
{\bf Proof.}
For any $\alpha\in C_{\rm 3-LieDer}^1(L;V)$, $\alpha$  is closed if and only if
$$\partial\alpha=({\rm d}\alpha,-\delta\alpha)=(0,0),$$
which is equivalent to that $\alpha[x,y,z]_L=\rho(x,y)\alpha(z)+\rho(y,z)\alpha(x)-\rho(x,z)\alpha(y)$ and $\alpha(\theta_L(x))=\theta_V(\alpha (x))$ for $x,y,z\in L$. Moreover, since there is no exact 1-cochain, we get the conclusion. \hfill $\Box$

\subsection{Maurer-Cartan characterizations of $3$-$\LD$ pairs }
\emptycomment{Let $\huaV$ be a graded vector space. The {\bf suspension operator $s$} assigns to $\huaV$ to the graded vector space $s\huaV$ with $(s\huaV)^i=\huaV^{i-1}$.
\begin{lemma}{\rm(\cite{TFS})}\label{lem:new graded Lie}
  Let $(\g,[\cdot,\cdot])$ be a graded Lie algebra. Define a linear map ${\rm ad}:\g\to{\rm gl}(s \g)$ by
  $${\rm ad}_x(su)=(-1)^{|x|}s[x,u]$$
  for homogenous elements $x,u\in\g$.
  Then $(\g;{\rm ad})$ is a representation of the graded Lie algebra $(\g,[\cdot,\cdot])$. Furthermore, we obtain a semi-direct product graded Lie algebra $\g\ltimes_{\rm ad} s\g$.
\end{lemma}}

A permutation $\sigma \in \mathbb S_n$ is called an $(i,n-i)$-shuffle if $\sigma(1) < \cdots <\sigma(i)$ and $\sigma(i+1) < \cdots < \sigma(n)$. If $i=0$ or $i=n$, we assume $\sigma={\rm Id}$. The set of all $(i,n-i)$-shuffles will be denoted by $\mathbb S_{(i,n-i)}$.

Let $L$ be a vector space. Denote by $C^\ast (L;L)=\oplus^\infty _{n=1}C^n(L;L)$, where $C^n(L;L)=\Hom(\otimes^{n-1}(\wedge^2 L)\wedge L,L)$. It was shown in \cite{Ro} that the graded vector space $C^\ast (L;L)$ equipped with the following bracket
\begin{eqnarray}\label{eq:graded bracket}
[P,Q]_{\rm 3Lie}=P\circ Q-(-1)^{pq}Q\circ P, \quad \forall P \in C^{p+1}(L;L),Q \in C^{q+1}(L;L),
\end{eqnarray}
is a degree $-1$ graded Lie algebra, where $P \circ Q \in C^{p+q+1}(L;L)$ is defined by
{\footnotesize
\begin{eqnarray}
&& \nonumber (P \circ Q)(\mathfrak{X}_1,...,\mathfrak{X}_{p+q},z)\\
~ \nonumber&=&\sum_{k=1}^{p}(-1)^{(k-1)q}\sum_{\sigma \in \mathbb{S}(k-1,q)}(-1)^{\sigma}P\big(\mathfrak{X}_{\sigma(1)},...,\mathfrak{X}_{\sigma(k-1)},Q(\mathfrak{X}_{\sigma(k)},...,\mathfrak{X}_{\sigma(k+q-1)},x_{k+q})\wedge y_{k+q},\mathfrak{X}_{k+q+1},...,\mathfrak{X}_{p+q},z\big)\\
~ \nonumber&&+\sum_{k=1}^{p}(-1)^{(k-1)q}\sum_{\sigma \in \mathbb{S}(k-1,q)}(-1)^{\sigma}P\big(\mathfrak{X}_{\sigma(1)},...,\mathfrak{X}_{\sigma(k-1)},x_{k+q}\wedge Q(\mathfrak{X}_{\sigma(k)},...,\mathfrak{X}_{\sigma(k+q-1)},y_{k+q}),\mathfrak{X}_{k+q+1},...,\mathfrak{X}_{p+q},z\big)\\
~ && +\sum_{\sigma \in \mathbb{S}(p,q)}(-1)^{pq}(-1)^{\sigma}P\big(\mathfrak{X}_{\sigma(1)},...,\mathfrak{X}_{\sigma(p)},Q(\mathfrak{X}_{\sigma(p+1)},...,\mathfrak{X}_{\sigma(p+q)},z)\big)
\end{eqnarray}}
for all $\mathfrak{X}_i=x_i \wedge y_i \in \wedge^2 L, ~1 \le i \le p+q, ~z \in L$.
In particular, $\pi :\wedge^3 L \to L$ defines a $3$-Lie algebra structure on $L$ if and only if $[\pi,\pi]_{\rm3Lie}=0$, i.e. $\pi$ is a Maurer-Cartan element of  graded Lie algebra $(C^{\bullet}(L;L),[\cdot,\cdot]_{\rm 3Lie})$. Moreover, the coboundary map $\dM$ of the $3$-Lie algebra with coefficients in the adjoint representation can be given by
$$\dM f_p=(-1)^{p-1}[\pi,f_p]_{\rm 3Lie}, \quad \forall f_p \in C^p(L;L).$$

\begin{prop}
 Let $L$ be a vector space. The bracket $[\cdot,\cdot]_{\rm 3-LieDer}:\mathcal{C}^{p}_{{\rm 3-LieDer}}(L;L)\times \mathcal{C}^{q}_{{\rm 3-LieDer}}(L;L)\to \mathcal{C}^{p+q-1}_{{\rm 3-LieDer}}(L;L)$ given by
  \begin{equation}
    [(\alpha_p,\beta_{p-1}),(\alpha_q,\beta_{q-1})]_{\rm 3-LieDer}=([\alpha_p,\alpha_q]_{\rm3Lie},(-1)^{p+1}[\alpha_p,\beta_{q-1}]_{\rm3Lie}+[\beta_{p-1},\alpha_q]_{\rm3Lie})
  \end{equation}
  defines a degree $-1$ graded Lie bracket on the graded vector space $\mathcal{C}^{\ast}_{{\rm 3-LieDer}}(L;L)$.  Moreover, its Maurer-Cartan elements are precisely the $3$-$\LD$ pairs on $L$.
\end{prop}
{\bf Proof.}The proof of $\mathcal{C}^{\ast}_{{\rm 3-LieDer}}(L;L)$ equipped with the bracket $[\cdot,\cdot]_{\rm 3-LieDer}$ is a degree $-1$ graded Lie algebra following as that  $(C^\ast (L;L),[\cdot,\cdot]_{\rm3Lie})$ is a degree $-1$ graded Lie algebra. For $\omega\in{\rm Hom}(\wedge^3L,L)$ and $\varphi\in{\rm Hom}(L,L)$, then we have
\begin{equation*}
   [(\omega,\varphi),(\omega,\varphi)]_{\rm 3-LieDer}=([\omega,\omega]_{\rm3Lie},-2[\omega,\varphi]_{\rm3Lie}).
\end{equation*}
It is obvious that $(\omega,\varphi)$ is a Maurer-Cartan element if and only if $[\omega,\omega]_{\rm3Lie}=0$ and $[\omega,\varphi]_{\rm3Lie}=0$. For $x,y,z\in L$, we have
\begin{eqnarray*}
  [\omega,\varphi]_{\rm3Lie}(x,y,z)=\omega(\varphi(x),y,z)+\omega(x,\varphi(y),z)+\omega(x,y,\varphi(z))-\varphi(\omega(x,y,z)).
\end{eqnarray*}
Thus $(\omega,\varphi)$ is a Maurer-Cartan element if and only if $(L,\omega,\varphi)$ is a $3$-$\LD$ pair.
 \hfill $\Box$
\vspace{3mm}

The following proposition shows that the coboundary operator $\partial$ with coefficient in the adjoint representation of a $3$-$\LD$ pair $(L,\omega,\varphi)$ can be presented by the degree $-1$ graded Lie algebra $(\mathcal{C}^{\ast}_{{\rm 3-LieDer}}(L;L),[\cdot,\cdot]_{\rm 3-LieDer})$.
\begin{prop}
  Let $(L,\omega,\varphi)$  be a $3$-$\LD$ pair. Then we have
  \begin{equation}
  \partial(\alpha_p,\beta_{p-1})=(-1)^{p-1}[(\omega,\varphi),(\alpha_p,\beta_{p-1})]_{\rm 3-LieDer},\quad \forall~(\alpha_p,\beta_{p-1})\in \mathcal{C}^{p}_{{\rm 3-LieDer}}(L;L).
  \end{equation}
\end{prop}
{\bf Proof.}
It follows by a direct calculation. We omit the details.
 \hfill $\Box$
\section{Deformations of {\rm 3-LieDer} pairs}\label{s3}
In this section, we study  one-parameter formal deformations of {\rm 3-LieDer} pairs, in which the $3$-Lie algebra and the distinguished derivation are simultaneous deformations.

\begin{definition}\label{1pfd}
 A {\bf one-parameter formal deformation} of a {\rm 3-LieDer} pair $(L,[\cdot,\cdot,\cdot]_L,\theta_L)$ is a pair of formal power series $(f_t,g_t)$ of the form
	\begin{align}
	f_t=\sum\limits_{i\geq0}f_it^i,\ \ g_t=\sum\limits_{i\geq0}g_it^i,\label{ftgtdef}
	\end{align}
	where $f_0=[\cdot,\cdot,\cdot]_L$ and $g_0=\theta_L$ such that $(L[[t]],f_t,g_t)$ is a {\rm 3-LieDer} pair over $\F[[t]]$.
\end{definition}

By a direct calculation, we have
\begin{lemma}\label{ftgtiff}
	A pair $(f_t,g_t)$ is a one-parameter formal deformation of $(L,\theta_L)$ if and only if the following equations hold:
	\begin{align}
	0=&\!\ \sum\limits_{\mbox{\tiny$\begin{array}{c}
			i+j=n\\
			n\geq0\end{array}$}} \Big(f_i(f_j(x_1,x_2,x_3),x_4,x_5)+f_i(x_3,f_j(x_1,x_2,x_4),x_5)\notag\\
	&\ \ +f_i(x_3,x_4,f_j(x_1,x_2,x_5))-f_i(x_1,x_2,f_j(x_3,x_4,x_5))\Big), \label{ftequ}\\
	0=&\!\ \sum\limits_{\mbox{\tiny$\begin{array}{c}
			i+j=n\\
			n\geq0\end{array}$}}\Big(g_i(f_j(x_1,x_2,x_3))-f_j(g_i(x_1),x_2,x_3)\notag\\
	&\ \ -f_j(x_1,g_i(x_2),x_3)-f_j(x_1,x_2,g_i(x_3))\Big).\label{gtequ}
	\end{align}
\end{lemma}

In particular,  for $n=1$,   (\ref{ftequ}) and (\ref{gtequ}) imply that
\begin{align}
0=&\!\ f_1(f_0(x_1,x_2,x_3),x_4,x_5)+f_1(x_3,f_0(x_1,x_2,x_4),x_5)+f_1(x_3,x_4,f_0(x_1,x_2,x_5))\notag\\
&\!\ -f_1(x_1,x_2,f_0(x_3,x_4,x_5))+f_0(f_1(x_1,x_2,x_3),x_4,x_5)+f_0(x_3,f_1(x_1,x_2,x_4),x_5)\notag\\
&\!\ +f_0(x_3,x_4,f_1(x_1,x_2,x_5))-f_0(x_1,x_2,f_1(x_3,x_4,x_5)),\label{n1f}\\
0=&\!\ g_1(f_0(x_1,x_2,x_3))-f_0(g_1(x_1),x_2,x_3)-f_0(x_1,g_1(x_2),x_3)\notag\\
&\!\ -f_0(x_1,x_2,g_1(x_3))+g_0(f_1(x_1,x_2,x_3))-f_1(g_0(x_1),x_2,x_3)\notag\\
&\!\ -f_1(x_1,g_0(x_2),x_3)-f_1(x_1,x_2,g_0(x_3)).\label{n1g}
\end{align}
Note that (\ref{n1f}) implies that ${\rm d}\!\!\ f_1=0$ and (\ref{n1g}) implies that ${\rm d}g_1+\delta f_1=0$. Thus we have

\begin{lemma}\label{f1g1z2}
$(f_1,g_1)$ is a $2$-cocycle, i.e.	$(f_1,g_1)\in\mathcal{Z}_{\rm 3-LieDer}^2(L;L)$.
\end{lemma}
{\bf Proof.}
Since $\partial(f_1,g_1)=({\rm d}f_1,{\rm d}g_1+\delta f_1)=(0,0)$, the conclusion follows. \hfill $\Box$
\vspace{3mm}

Similar to the proof of Lemma \ref{f1g1z2}, we can show that
\begin{prop}\label{fngnz2}
	If $(f_i,g_i)=0$, $1\leq i<n$, then $(f_n,g_n)\in\mathcal{Z}_{\rm 3-LieDer}^2(L;L)$.
\end{prop}

A $2$-cochain $(f_n,g_n)$ is called the {\bf $n$-infinitesimal} of $(f_t,g_t)$ if $(f_i,g_i)=0$ for all $1\leq i<n$. In particular, the $2$-cocycle $(f_1,g_1)$ is called the {\bf infinitesimal} (or $1$-infinitesimal) of $(f_t,g_t)$.

Therefore, by Proposition \ref{fngnz2}, we see that the $n$-infinitesimal $(f_n,g_n)$ of $(f_t,g_t)$ is a $2$-cocycle whenever $n\geq1$.

\begin{definition}\label{ftftpequ}
	Suppose that $(f_t,g_t)$ and $(f_{t}',g_{t}')$ are one-parameter formal deformations of a {\rm 3-LieDer} pair $(L,\theta_L)$. $(f_t,g_t)$ and $(f_{t}',g_{t}')$ are called {\bf equivalent} (i.e., $(f_t,g_t)\sim(f_{t}',g_{t}')$), if there exists a linear isomorphism $\Phi_t=\sum_{i\geq0}\phi_it^i:(L[[t]],f_t,g_t)\rightarrow (L[[t]],f_t,g_t)$, where $\phi_i\in C^1(L;V)$ and $\phi_0=id_L$, satisfying that
	\begin{align}
	\Phi_t\circ f_{t}'=f_t\circ (\Phi_t\times\Phi_t\times\Phi_t),\ \ \Phi_t\circ g_{t}'=g_t\circ\Phi_t.\label{pairequ}
	\end{align}
\end{definition}

Motivated by Proposition \ref{fngnz2}, we have the following
\begin{prop}\label{equsame}
	The $n$-infinitesimals of two equivalent one-parameter formal deformations of a {\rm 3-LieDer} pair $(L,\theta_L)$ belong to the same cohomology class.
\end{prop}
{\bf Proof.}
Let $(f_t,g_t)$ and $(f_{t}',g_{t}')$ be one-parameter formal deformations of $(L,\theta_L)$. It suffices to show that $[(f_n,g_n)]=[(f_{n}',g_{n}')]\in\mathcal{H}_{\rm 3-LieDer}^2(L;L)$. By exacting the coefficients of $t^n$, we see that (\ref{pairequ}) is equivalent to
\begin{align*}
\sum\limits_{\mbox{\tiny$\begin{array}{c}
		i+j=n\\
		n\geq0\end{array}$}}\phi_i(f_j(x,y,z))&=\sum\limits_{\mbox{\tiny$\begin{array}{c}
		i+j+k+l=n\\
		n\geq0\end{array}$}}f_{l}'(\phi_i(x),\phi_j(y),\phi_k(z)),\\
\sum\limits_{\mbox{\tiny$\begin{array}{c}
		i+j=n\\
		n\geq0\end{array}$}}\phi_i(g_j(x))&=\sum\limits_{\mbox{\tiny$\begin{array}{c}
		i+j=n\\
		n\geq0\end{array}$}}g_{i}'(\phi_j(x))
\end{align*}
for all $x,y,z\in L$.

Since $\phi_o=id_L$, we have
\begin{align*}
f_n(x,y,z)-f_{n}'(x,y,z)=&\!\ -\phi_1(f_0(x,y,z))+f_n(\phi_1(x),y,z)+f_n(x,\phi_1(y),z)\\
&\!\ +f_n(x,y,\phi_1(z)),\\
g_n(x)-g_{n}'(x)=&\!\ -\phi_1(\theta_L(x))+\theta_L(\phi_1(x)),
\end{align*}
which implies that $f_n-f_{n}'={\rm d}\phi_1$ and $g_n-g_{n}'=-\delta\phi_1$. Thus we have
\begin{align*}
(f_n,g_n)-(f_{n}',g_{n}')=({\rm d}\phi_1,-\delta\phi_1)=\partial\phi_1.
\end{align*}
This means that $[(f_n,g_n)]=[(f_{n}',g_{n}')]\in\mathcal{H}_{\rm 3-LieDer}^2(L;L)$. \hfill $\Box$
\vspace{3mm}

Next we shall investigate the rigidity of {\rm 3-LieDer} pairs.
\begin{definition}
	Let $(f_{t},g_{t})$ be a one-parameter formal deformation of a {\rm 3-LieDer} pair $(L,\theta_L)$. If $(f_t,g_t)\sim(f_0,g_0)$, then $(f_t,g_t)$ is said to be {\bf trivial}. Moreover, $(L,\theta_L)$ is called {\bf rigid} if every one-parameter formal deformation of $(L,\theta_L)$ is trivial.
\end{definition}

\begin{theorem}
	Let $(L,\theta_L)$ be a {\rm 3-LieDer} pair. If $\mathcal{H}_{\rm 3-LieDer}^2(L;L)=0$, then $(L,\theta_L)$ is rigid.
\end{theorem}
{\bf Proof.}
Let $(f_t,g_t)$ be a one-parameter formal deformation of $(L,\theta_L)$. Set $f_t=f_0+\sum\limits_{i\geq r}f_it^i,\ \ g_t=g_0+\sum\limits_{i\geq r}g_it^i$. By Proposition \ref{fngnz2}, $(f_r,g_r)\in\mathcal{Z}_{\rm 3-LieDer}^2(L;L)$. Since $\mathcal{H}_{\rm 3-LieDer}^2(L;L)=0$, there is a map $\phi_r\in C^1(L;L)$ such that $(f_r,g_r)=-\partial\phi_r=(-{\rm d}\phi_r,\delta\phi_r)$. Then $\Phi_t:=id_L+\phi_rt^r$ is a linear isomorphism of $L$. Define $f_{t}',g_{t}'$ as
\begin{align}
f_{t}'=\Phi_t^{-1}\circ f_t\circ(\Phi_t\times\Phi_t\times\Phi_t),\ \ \ \ g_{t}'=\Phi_t^{-1}\circ g_t\circ\Phi_t. \label{ftpdef}
\end{align}
It is straightforward to check that $(f_{t}',g_{t}')$ is a one-parameter formal deformation of $(L,\theta_L)$. Thus, by (\ref{ftpdef}) and Definition \ref{ftftpequ}, one has $(f_{t}',g_{t}')\sim(f_t,g_t)$.

 By a direct calculation, we have
\begin{align}
f_{t}'(x,y,z)&=f_0(x,y,z)+\Big(f_r(x,y,z)-\phi_r(f_0(x,y,z))+f_r(\phi_r(x),y,z)\notag\\
&\ \ \ \ +f_r(x,\phi_r(y),z)+f_r(x,y,\phi_r(z))\Big)t^r+\cdots,\label{ftpr}\\
g_{t}'(x)&=g_0(x)+\Big(g_r(x)+g_0(\phi_r(x))-\phi_r(g_0(x))\Big)t^r+\cdots.\label{gtpr}
\end{align}
Since $f_r=-\partial\phi_r$ and $g_r=\delta\phi_r$, (\ref{ftpr}) and (\ref{gtpr}) reduce to
\begin{align*}
f_{t}'(x,y,z)&=f_0(x,y,z)+\sum\limits_{i\geq{r+1}}f_{i}'(x,y,z)t^{i},\ g_{t}'(x)=g_0(x)+\sum\limits_{i\geq{r+1}}g_{i}'(x)t^{i}.
\end{align*}
Then by repeating the argument, we obtain that $(f_t,g_t)\sim(f_0,g_0)$. \hfill $\Box$

\begin{definition}\label{deforn}
	A {\bf  deformation of order $n$} of a {\rm 3-LieDer} pair $(L,\theta_L)$ is a pair $(f_t,g_t)$ such that $f_t=\sum_{i=0}^{n}f_it^i$ and $\  g_t=\sum_{i=0}^{n}g_it^i$ endow the $\F[[t]]/(t^{n+1})$-module $L[[t]]/(t^{n+1})$ the  {\rm 3-LieDer} pair structure with
 $(f_0,g_0)=([\cdot,\cdot,\cdot]_L,\theta_L)$. Furthermore,  if there exists a 2-cochain $(f_{n+1},g_{n+1})$ belongs to $\mathcal{C}_{\rm 3-LieDer}^2(L;L)$ such that $(f_{t}',g_{t}')$ defined by
\begin{align}
f_{t}'=f_t+f_{n+1},\ g_{t}'=g_t+g_{n+1}
\end{align}
is  a deformation of order $n+1$  of $(L,\theta_L)$,  we say that $(f_t,g_t)$ is {\bf extensible}.
\end{definition}

Let $(f_t,g_t)$ be a deformation of order $n$ of a {\rm 3-LieDer} pair $(L,\theta_L)$. For any $X_1=x_1\wedge x_2, X_2=x_3\wedge x_4\in\wedge^2L$ and $x_5\in L$, define a $3$-cochain $(\Omega_{n+1}^3,\Omega_{n+1}^2)\in\mathcal{C}_{\rm 3-LieDer}^3(L;L)$ as
\begin{align}
\Omega_{n+1}^3(X_1,X_2,x_5)=&\sum\limits_{\mbox{\tiny$\begin{array}{c}
		i+j=n+1\\
		i,j>0\end{array}$}}\!\!\Big(f_i(f_j(x_1,x_2,x_3),x_4,x_5)+f_i(x_3,f_j(x_1,x_2,x_4),x_5)\notag\\
&\ \ +f_i(x_3,x_4,f_j(x_1,x_2,x_5))-f_i(x_1,x_2,f_j(x_3,x_4,x_5))\Big), \notag\\
\Omega_{n+1}^2(X_1,x_3)=& \sum\limits_{\mbox{\tiny$\begin{array}{c}
		i+j=n+1\\
		i,j>0\end{array}$}}\Big(g_i(f_j(x_1,x_2,x_3))-f_j(g_i(x_1),x_2,x_3)\notag\\
&\ \ -f_j(x_1,g_i(x_2),x_3)-f_j(x_1,x_2,g_i(x_3))\Big).\notag
\end{align}

\emptycomment{By the  graded Lie bracket on $C^\ast (L;L)$ given by  (\ref{eq:graded bracket}), $\Omega_{n+1}^3$ and $\Omega_{n+1}^2$ can be rewritten as
\begin{align}
\Omega_{n+1}^3&=\frac{1}{2}\sum\limits_{\mbox{\tiny$\begin{array}{c}
		i+j=n+1\\
		i,j>0\end{array}$}}[f_i,f_j]_{\rm 3Lie},\label{O3}\\
\Omega_{n+1}^2&=\sum\limits_{\mbox{\tiny$\begin{array}{c}
		i+j=n+1\\
		i,j>0\end{array}$}}[g_i,f_j]_{\rm 3Lie}.\label{O2}
\end{align}
}

\emptycomment{\begin{lemma}\label{O3O2sim}
	(1) $\Omega_{n+1}^3=-[f_0,f_{n+1}]^{3Lie}$, or equivalently, $\Omega_{n+1}^3={\rm d}f_{n+1}$.\\
	(2) $\Omega_{n+1}^2={\rm d}g_{n+1}+\delta f_{n+1}$.
\end{lemma}
{\bf Proof.}
By (\ref{O3}) and (\ref{ftequ}), for any $x_i\in L$ we get
\begin{align*}
&\sum\limits_{\mbox{\tiny$\begin{array}{c}
		i+j=n\\
		n>0\end{array}$}} \Big(f_i(f_j(x_1,x_2,x_3),x_4,x_5)+f_i(x_3,f_j(x_1,x_2,x_4),x_5)\\
&\ \ +f_i(x_3,x_4,f_j(x_1,x_2,x_5))-f_i(x_1,x_2,f_j(x_3,x_4,x_5))\Big)\\
=&-\Big(f_0(f_{n+1}(x_1,x_2,x_3),x_4,x_5)+f_0(x_3,f_{n+1}(x_1,x_2,x_4),x_5)\\
&\ \ +f_0(x_3,x_4,f_{n+1}(x_1,x_2,x_5))-f_0(x_1,x_2,f_{n+1}(x_3,x_4,x_5))\\
&\ \ +f_{n+1}(f_0(x_1,x_2,x_3),x_4,x_5)+f_{n+1}(x_3,f_0(x_1,x_2,x_4),x_5)\\
&\ \ +f_{n+1}(x_3,x_4,f_0(x_1,x_2,x_5))-f_{n+1}(x_1,x_2,f_0(x_3,x_4,x_5))\Big),
\end{align*}
which is equivalent to $\Omega_{n+1}^3=-[f_0,f_{n+1}]^{3Lie}$ due to (\ref{NRbr}).

Similarly, we can show that
\begin{align*}
\Omega_{n+1}^2&=[f_0,g_{n+1}]^{3Lie}-[g_0,f_{n+1}]^{3Lie},
\end{align*}
which implies that $\Omega_{n+1}^2=dg_{n+1}+\delta f_{n+1}$ as required. \hfill $\Box$}
\begin{prop}
	The $3$-cochain $(\Omega_{n+1}^3,\Omega_{n+1}^2)$ is closed, i.e. $\partial (\Omega_{n+1}^3,\Omega_{n+1}^2)=0.$
\end{prop}
{\bf Proof.}
It is straightforward to check that $\Omega_{n+1}^3={\rm d}f_{n+1}$ and $\Omega_{n+1}^2={\rm d}g_{n+1}+\delta f_{n+1}$. It is obvious that
${\rm d}\Omega_{n+1}^3={\rm d}({\rm d}f_{n+1})=0$. By Lemma \ref{dd=dd}, we also have
\begin{align*}
{\rm d}\Omega_{n+1}^2-\delta\Omega_{n+1}^3={\rm d}({\rm d}g_{n+1}+\delta f_{n+1})-\delta({\rm d}f_{n+1})=({\rm d}\circ\delta-\delta\circ {\rm d})f_{n+1}=0.
\end{align*}
Thus $\partial(\Omega_{n+1}^3,\Omega_{n+1}^2)=({\rm d}\Omega_{n+1}^3,{\rm d}\Omega_{n+1}^2-\delta\Omega_{n+1}^3)=(0,0)$. The proof is completed. \hfill $\Box$
\vspace{3mm}

The main result of this section is given in the following
\begin{theorem}\label{ntri}
	Let $(f_t,g_t)$ be a deformation of order $n$ of a {\rm 3-LieDer} pair $(L,\theta_L)$. Then $(f_t,g_t)$ is extensible if and only if the cohomology class $[(\Omega_{n+1}^3,\Omega_{n+1}^2)]$ in $\mathcal{H}_{\rm 3-LieDer}^3(L;V)$ is trivial.
\end{theorem}
{\bf Proof.}
Suppose that $(f_t,g_t)$ is extensible, i.e., $(f_{t}',g_{t}')$ is a deformation of order $n+1$, where
\begin{align*}
f_{t}'=f_t+f_{n+1},\ g_{t}'=g_t+g_{n+1}.
\end{align*}
It is straightforward to check that $\Omega_{n+1}^3={\rm d}f_{n+1}, \Omega_{n+1}^2={\rm d}g_{n+1}+\delta f_{n+1}$, which imply that
\begin{align*}
(\Omega_{n+1}^3,\Omega_{n+1}^2)=({\rm d}f_{n+1},{\rm d}g_{n+1}+\delta f_{n+1})=\partial(f_{n+1},g_{n+1}).
\end{align*}
Thus the cohomology class $[(\Omega_{n+1}^3,\Omega_{n+1}^2)]$ is trivial.

Conversely, suppose that $(\Omega_{n+1}^3,\Omega_{n+1}^2)$ is trivial. Then there exists a 2-cochain $(f_{n+1},g_{n+1})\in\mathcal{C}_{\rm 3-LieDer}^2(L;L)$ such that $(\Omega_{n+1}^3,\Omega_{n+1}^2)=\partial(f_{n+1},g_{n+1})$. Set $f_{t}'=f_t+f_{n+1},\ g_{t}'=g_t+g_{n+1}$, then by a direct check we see that $(f_{t}',g_{t}')$ is a {3-LieDer} pair of $(L,\theta_L)$ and hence $(f_{t}',g_{t}')$ is a deformation of order $n+1$ due to Definition \ref{deforn}. Namely, $(f_t,g_t)$ is extensible. The proof is finished. \hfill $\Box$

\section{Abelian extensions of {\rm 3-LieDer} pairs}\label{abel}
In this section, we study abelian extensions of {\rm 3-LieDer} pairs and show that equivalent abelian extensions of {\rm 3-LieDer} pairs can be classified by the second cohomology groups.

\begin{definition}\label{aep}
	Let $(A,\theta_A)$ and $(B,\theta_B)$ be two {\rm 3-LieDer} pairs. An {\bf abelian extension} of $(B,\theta_B)$ by $(A,\theta_A)$ is an exact sequence of {\rm 3-LieDer} pair morphisms
	\begin{eqnarray}
	\begin{CD}
	0 @>>> {A} @>{i}>> L @>\pi>> B @>>> 0 \\
	&& @V  \theta_A VV  @V  \theta_L VV  @V  \theta_B VV   \\
	0 @>>> A @>{i}>> L @>\pi>> B @>>> 0
	\end{CD}
	\end{eqnarray}
	such that $A$ is an abelian ideal of $L$, i.e., $[u,v,\cdot]_L=0$, $\forall\ u,v\in A$. Denote this abelian extension by $\varepsilon_{(L,\theta_L)}$.  A {\bf section} of an abelian extension of $(B,\theta_B)$ by $(A,\theta_A)$  is a linear map $s:B\rightarrow L$ such that
$\pi\circ s=id_B$.
\end{definition}

Let $(B,\theta_B)$ be a {\rm 3-LieDer} pair and $(A;\rho,\theta_A)$ a representation of $(B,\theta_B)$. Suppose that $(\psi,\lambda)\in\mathcal{C}_{\rm 3-LieDer}^2(B;A)$. Define $[\cdot,\cdot,\cdot]_{(\rho,\psi)}:\wedge^3(B\oplus A)\rightarrow B\oplus A$ and $\theta_{\lambda}:B\oplus A\rightarrow B\oplus A$ as
\begin{align}
[x_1+a_1,x_2+a_2,x_3+a_3]_{(\rho,\psi)}=&\!\ [x_1,x_2,x_3]_B+\psi(x_1,x_2,x_3)+\rho(x_1,x_2)(a_3)\notag\\
&\!\ +\rho(x_2,x_3)(a_1)+\rho(x_3,x_1)(a_2) \label{semipro},\\
\theta_{\lambda}(x+a)=&\!\ \theta_B(x)+\lambda(x)+\theta_A(a)\label{thlam}
\end{align}
for all  $x,x_i\in B, a,a_i\in A.$  Denote $B\oplus A$ with the bracket operation \eqref{semipro} by $B\ltimes_{\rho,\psi} A$.
\begin{prop}\label{boa}
	With the above notation. Then $(B\ltimes_{\rho,\psi} A,\theta_{\lambda})$ is a {\rm 3-LieDer} pair if and only if $(\psi,\lambda)\in\mathcal{Z}_{\rm 3-LieDer}^2(B;A)$. In this case,
	$0\rightarrow (A,\theta_A)\hookrightarrow (B\ltimes_{\rho,\psi} A,\theta_{\lambda})\overset{p}{\rightarrow} (B,\theta_B)\rightarrow 0$ is an abelian extension, where $p$ is the canonical projection.
\end{prop}
{\bf Proof.}
 Suppose that $(B\oplus A,\theta_{\lambda})$ is a {\rm 3-LieDer} pair. It is straightforward to check that $B\oplus A$ equipped with the bracket $[\cdot,\cdot,\cdot]_{(\rho,\psi)}$ is a 3-Lie algebra if and only if ${\rm d}\psi=0$. Moreover, by the fact that $\theta_\lambda$ is a derivation on $B\ltimes_{\rho,\psi} A$, $\theta_B$ is a derivation on $B$ and $(A;\rho,\theta_A)$ is a representation of $(B,\theta_B)$, for any $x_i\in B, a_i\in A$, we have
\begin{eqnarray*}
&&[\theta_\lambda(x_1+a_1),x_2+a_2,x_3+a_3]_{(\rho,\psi)}+[x_1+a_1,\theta_\lambda(x_2+a_2),x_3+a_3]_{(\rho,\psi)}\\
&& +[x_1+a_1,x_2+a_2,\theta_\lambda(x_3+a_3)]_{(\rho,\psi)}-\theta_\lambda([x_1+a_1,x_2+a_2,x_3+a_3]_{(\rho,\psi)})\\
&=&\mathop{\circlearrowleft}\limits_{x_1,x_2,x_3}[\theta_B(x_1),x_2,x_3]_B-\theta_B([x_1,x_2,x_3]_B)-\lambda([x_1,x_2,x_3]_B)\\
&&+\mathop{\circlearrowleft}\limits_{x_1,x_2,x_3}\rho(x_1,x_2)\lambda(x_3)+\mathop{\circlearrowleft}\limits_{x_1,x_2,x_3}\psi(\theta_B(x_1),x_2,x_3)-\theta_A(\psi(x_1,x_2,x_3))\\
&& +\rho(\theta_B(x_1),x_2)(a_3)+\rho(x_1,\theta_B(x_2))(a_3)+\rho(x_1,x_2)\theta_A(a_3)-\theta_A(\rho(x_1,x_2)(a_3))\\
&& +\rho(\theta_B(x_2),x_3)(a_1)+\rho(x_2,\theta_B(x_3))(a_1)+\rho(x_2,x_3)\theta_A(a_1)-\theta_A(\rho(x_2,x_3)(a_1))\\
&& +\rho(\theta_B(x_3),x_1)(a_2)+\rho(x_3,\theta_B(x_1))(a_2)+\rho(x_3,x_1)\theta_A(a_2)-\theta_A(\rho(x_3,x_1)(a_2))\\
&=&-\lambda([x_1,x_2,x_3]_B)+\mathop{\circlearrowleft}\limits_{x_1,x_2,x_3}\rho(x_1,x_2)\lambda(x_3) +\mathop{\circlearrowleft}\limits_{x_1,x_2,x_3}\psi(\theta_B(x_1),x_2,x_3)-\theta_A(\psi(x_1,x_2,x_3))\\
&=&({\rm d}\lambda+\delta\psi)(x_1,x_2,x_3)=0,
\end{eqnarray*}
where $\mathop{\circlearrowleft}\limits_{x_1,x_2,x_3}$ denotes the summation over the cyclic permutations of ${x_1},{x_2},{x_3}$. Furthermore, we have
\begin{align*}
\partial(\psi,\lambda)=({\rm d}\psi,{\rm d}\lambda+\delta\psi)=(0,0),
\end{align*}
which implies that $(\psi,\lambda)\in\mathcal{Z}_{\rm 3-LieDer}^2(B;A)$.

The converse can be proved similarly. We omit the details. \hfill $\Box$
\vspace{3mm}

Let $\varepsilon_{(L,\theta_L)}$ be an abelian extension of $(B,\theta_B)$ by $(A,\theta_A)$ and $s:B\rightarrow L$ a section. Define $\rho:\wedge^2B\rightarrow {\rm End}(A)$, $\omega:\wedge^3B\rightarrow A$ and $\mu:B\rightarrow A$ by
\begin{align}
\rho(x,y)(a)&=[s(x),s(y),a]_L,\label{rhoext}\\
\omega(x,y,z)&=[s(x),s(y),s(z)]_L-s([x,y,z]_B),\ \ \label{omega}\\
\mu(x)&=\theta_L(s(x))-s(\theta_B(x)),\ \forall\ \ x,y,z\in B,a\in A.\label{mu}
\end{align}
Note that $\rho$ does not depend on the choice of sections.
\begin{prop}\label{extabel}
	Let $\varepsilon_{(L,\theta_L)}$ be an abelian extension of $(B,\theta_B)$ by $(A,\theta_A)$. Then $(A;\rho,\theta_A)$ is a representation of the {\rm 3-LieDer} pair $(B,\theta_B)$ and $(\omega,\mu)\in\mathcal{Z}_{\rm 3-LieDer}^2(B;A)$. Moreover, the cohomology class $[(\omega,\mu)]\in\mathcal{H}_{\rm 3-LieDer}^2(B;A)$ does not depend on the choice of sections.
\end{prop}
{\bf Proof.}
It is routine to check that  $\rho$ is a representation of the $3$-Lie algebra $B$ on $
A$. We only need to show the following equality
\begin{align}
\theta_A(\rho(x,y)(a))=\rho(\theta_B(x),y)(a)+\rho(x,\theta_B(y))(a)+\rho(x,y)(\theta_A(a)),\quad\forall~ x,y\in B,a\in A.\label{rhoppf}
\end{align}
By a direct calculation, we have
\begin{eqnarray*}
&&\rho(\theta_B(x),y)(a)+\rho(x,\theta_B(y))(a)+\rho(x,y)(\theta_A(a))\\
&=&[s(\theta_B(x)),s(y),a]_L+[s(x),s(\theta_B(y)),a]_L+[s(x),s(y),\theta_L(a)]_L\\
&=&[\theta_L(s(x))-\mu(x),s(y),a]_L+[s(x),\theta_L(s(y))-\mu(y),a]_L+[s(x),s(y),\theta_L(a)]_L\\
&=&[\theta_L(s(x)),s(y),a]_L+[s(x),\theta_L(s(y)),a]_L+[s(x),s(y),\theta_L(a)]_L\\
&=&\theta_L([s(x),s(y),a]_L)\\
&=&\theta_A(\rho(x,y)(a)).
\end{eqnarray*}
Thus $(A;\rho,\theta_A)$ is a representation of the {\rm 3-LieDer} pair $(B,\theta_B)$.

Since $\varepsilon_{(L,\theta_L)}$ is an abelian extension of $(B,\theta_B)$ by $(A,\theta_A)$, by Proposition \ref{boa}, $(\omega,\mu)\in\mathcal{Z}_{\rm 3-LieDer}^2(B;A)$.  Let $s_1$ and $s_2$ be two sections of $\pi$. Set $\lambda(x)=s_1(x)-s_2(x)$, $x\in B$. Note that $\lambda\in {\rm Hom}(B,A)$. By a direct calculation, we obtain that $\omega_1-\omega_2={\rm d}\lambda$. Moreover,
\begin{align*}
\mu_1(x)-\mu_2(x)=&\!\ \theta_L(s_1(x))-s_1(\theta_B(x))-(\theta_L(s_2(x))-s_2(\theta_B(x)))\\
=&\!\ \theta_L(s_1(x)-s_2(x))-(s_1(\theta_B(x))-s_2(\theta_B(x)))\\
=&\!\ \theta_A(\lambda(x))-\lambda(\theta_B(x))\\
=&\!\ -(\delta\lambda)(x).
\end{align*}
Then we have $(\omega_1,\mu_1)-(\omega_2,\mu_2)=({\rm d}\lambda,-\delta\lambda)=\partial \lambda$. Thus the cohomology class $(\omega_1,\mu_1)$ and $(\omega_2,\mu_2)$ are in the same cohomological class.\hfill $\Box$
 \vspace{3mm}

\begin{definition}
	Two abelian extensions $\varepsilon_{(L_1,\theta_{L_1})}, \varepsilon_{(L_2,\theta_{L_2})}$ of $(B,\theta_B)$ by $(A,\theta_A)$ are called {\bf equivalent}, if there exists a {\rm 3-LieDer} pair morphism $\eta:(L_1,\theta_{L_1})\rightarrow(L_2,\theta_{L_2})$ such that we have the following commutative diagram:
	\begin{eqnarray}\label{aeequdia}
	\CD
	0 @>>> (A,\theta_A) @>{\imath_1}>> (L_1,\theta_{L_1}) @>\pi_1>> (B,\theta_B) @>>> 0 \\
	&&  @|   @V \eta VV @|      \\
	0 @>>> (A,\theta_A) @>{\imath_2}>> (L_2,\theta_{L_2}) @>\pi_2>> (B,\theta_B) @>>> 0.
	\endCD
	\end{eqnarray}
\end{definition}

Assume that $\varepsilon_{(L_1,\theta_{L_1})}$ and $\varepsilon_{(L_2,\theta_{L_2})}$ are equivalent abelian extensions of $(B,\theta_B)$ by $(A,\theta_A)$. It is not hard to see that these equivalent abelian extensions of {\rm 3-LieDer} pairs give the same representation $\rho$ defined by (\ref{rhoext}) of the $3$-Lie algebra $B$ on $A$. Thus $(A;\rho,\theta_A)$ is also a representation of the {\rm 3-LieDer} pair $(B,\theta_B)$, which does not depend on equivalence classes of abelian extensions of $(B,\theta_B)$ by $(A,\theta_A)$. Furthermore, we have

\begin{theorem}\label{equcla}
There is a one-to-one correspondence between equivalence classes of abelian extensions of $(B,\theta_B)$ by $(A,\theta_A)$ and the second cohomology group $\mathcal{H}_{\rm 3-LieDer}^2(B;A)$.
\end{theorem}
{\bf Proof.}
Let $\varepsilon_{(L_1,\theta_{L_1})}, \varepsilon_{(L_2,\theta_{L_2})}$ be two equivalent abelian extensions of $(B,\theta_B)$ by $(A,\theta_A)$. Choose a section $s_1:B\rightarrow L_1$ of $\pi_1$. Since
\begin{align*}
\pi_2\circ(\eta\circ s_1)=(\pi_2\circ\eta)\circ s_1=\pi_1\circ s_1=id_B,
\end{align*}
we obtain that $s_2:=\eta\circ s_1$ is a section of $\pi_2$. Let $\omega_i,\mu_i$ respectively given by (\ref{omega}),(\ref{mu}) corresponding to $s_i$, $i=1,2$. As $\eta$ is a {\rm 3-LieDer} pair morphism and $\eta|_A=id_A$, for any $x,y,z\in B$, we have
\begin{align*}
\omega_2(x,y,z)=&\!\ [s_2(x),s_2(y),s_2(z)]_{L_2}-s_2([x,y,z]_B)\\
=&\!\ [\eta(s_1(x)),\eta(s_1(x)),\eta(s_1(x))]_{L_2}-\eta(s_1([x,y,z]_B))\\
=&\!\ \eta([s_1(x),s_1(y),s_1(z)]_{L_1}-s_1([x,y,z]_B))\\
=&\!\ \eta(\omega_1(x,y,z))=\omega_1(x,y,z).
\end{align*}
Similarly, by $\eta\circ\theta_{L_1}=\theta_{L_2}\circ\eta$, we have
\begin{align*}
\mu_2(x)=&\!\ \theta_{L_2}(s_2(x))-s_2(\theta_{L_1}(x))\\
=&\!\ \theta_{L_2}(\eta(s_1(x)))-\eta(s_1(\theta_{L_1}(x)))\\
=&\!\ \eta(\theta_{L_1}(s_1(x)))-\eta(s_1(\theta_{L_1}(x)))\\
=&\!\ \eta(\theta_{L_1}(s_1(x))-s_1(\theta_{L_1}(x)))\\
=&\!\ \eta(\mu_1(x))=\mu_1(x).
\end{align*}
Hence we have $(\omega_1,\mu_1)=(\omega_2,\mu_2)$. By Proposition \ref{extabel}, the cohomology class $[(\omega_1,\mu_1)]$ of the abelian extension  $\varepsilon_{(L_1,\theta_{L_1})}$ and the cohomology class $[(\omega_2,\mu_2)]$ of the abelian extension  $\varepsilon_{(L_2,\theta_{L_2})}$ do not depend on the choice of sections. Thus equivalence abelian extensions of $(B,\theta_B)$ by $(A,\theta_A)$ give the same element in $\mathcal{H}_{\rm 3-LieDer}^2(B;A)$.

Conversely, let $[(\omega_1,\mu_1)]=[(\omega_2,\mu_2)]\in\mathcal{H}_{\rm 3-LieDer}^2(B;A)$. Then there exists a map $\lambda\in{\rm Hom}(B,A)$ such that
$$(\omega_1,\mu_1)-(\omega_2,\mu_2)=\partial\lambda=({\rm d}\lambda,-\delta\lambda).$$
By Proposition \ref{boa}, $0\rightarrow (A,\theta_A)\hookrightarrow (B\ltimes_{\rho,\omega_1} A,\theta_{\mu_1})\overset{p_1}{\rightarrow} (B,\theta_B)\rightarrow 0$ and $0\rightarrow (A,\theta_A)\hookrightarrow (B\ltimes_{\rho,\omega_2} A,\theta_{\mu_2})\overset{p_2}{\rightarrow} (B,\theta_B)\rightarrow 0$ are abelian extensions. Define $\eta:B\ltimes_{\rho,\omega_1} A \rightarrow B\ltimes_{\rho,\omega_2} A$ by
\begin{align}
\eta(x+a)=x+\lambda(x)+a,\ \forall\ x\in B,\ a\in A.
\end{align}
By a direct computation, we can show that $\eta$ is a {\rm 3-LieDer} pair morphism such that (\ref{aeequdia}) commutes. Thus these two abelian extensions are equivalent.  \hfill $\Box$

\section{Classification of skeletal and strict  {\rm $3$-Lie$2$Der} pairs}\label{s5}
In this section, first we recall the concept of 2-term 3-{$Lie_{\infty}$}-algebra, which is equivalent to $3$-Lie $2$-algebra. Then we introduce the definition of a {\rm $3$-Lie2Der} pair, which consists of a $3$-Lie $2$-algebra and a $2$-derivation. We show that skeletal {\rm $3$-Lie$2$Der} pairs can be classified by the triple which consists of a {\rm $3$-LieDer} pair, a representation and a $3$-cocycle. Analogous to crossed modules of $3$-Lie algebras (\cite{ZLS}), we introduce the definition of crossed modules of {\rm $3$-LieDer} pairs.  We show that there exists a one-to-one correspondence between strict {\rm $3$-Lie$2$Der} pairs and crossed modules of {\rm $3$-LieDer} pairs.

\begin{definition}\label{3lie2alg}{\rm(\cite{ZLS})}
	A {\bf $3$-Lie $2$-algebra}
 $\mathcal{V}=(V_1,V_0,d,l_3,l_5)$ consists of the following data:
	\begin{itemize}
\item[$\bullet$] a complex of vector spaces $V_1\overset{d}{\rightarrow} V_0$.
\item[$\bullet$] completely skew-symmetric trilinear maps $l_3:V_i\times V_j\times V_k\rightarrow V_{i+j+k}$, where $0\leq i+j+k\leq 1$.
\item[$\bullet$] a multilinear map $l_5:(\wedge^2 V_0)\otimes(\wedge^3 V_0)\rightarrow V_1$,
 \end{itemize}
	such that for any $x,y,x_i\in V_0$ and $u,v,w\in V_1$, the following equalities are satisfied:
	\begin{enumerate}
		\item[\rm(a)] $dl_3(x,y,u)=l_3(x,y,du)$,
		\item[\rm(b)] $l_3(u,v,w)=0;\ l_3(u,v,x)=0$,
		\item[\rm(c)] $l_3(du,v,x)=l_3(u,dv,x)$,
		\item[\rm(d)] $dl_5(x_1,x_2,x_3,x_4,x_5)=l_3(l_3(x_1,x_2,x_3),x_4,x_5)+l_3(x_3,l_3(x_1,x_2,x_4),x_5)$\\
		$+l_3(x_3,x_4,l_3(x_1,x_2,x_5))-l_3(x_1,x_2,l_3(x_3,x_4,x_5))$,
		\item[\rm(e)] $l_5(du,x_2,x_3,x_4,x_5)=l_3(l_3(u,x_2,x_3),x_4,x_5)+l_3(x_3,l_3(u,x_2,x_4),x_5)$\\
		$+l_3(x_3,x_4,l_3(u,x_2,x_5))-l_3(u,x_2,l_3(x_3,x_4,x_5))$,
		\item[\rm(f)] $l_5(x_1,x_2,du,x_4,x_5)=l_3(l_3(x_1,x_2,u),x_4,x_5)+l_3(u,l_3(x_1,x_2,x_4),x_5)$\\
		$+l_3(u,x_4,l_3(x_1,x_2,x_5))-l_3(x_1,x_2,l_3(u,x_4,x_5))$,
		\item[\rm(g)] $l_3(l_5(x_1,x_2,x_3,x_4,x_5),x_6,x_7)+l_3(x_5,l_5(x_1,x_2,x_3,x_4,x_6),x_7)+l_3(x_1,x_2,l_5(x_3,x_4,x_5,x_6,x_7))$\\
			$+l_3(x_5,x_6,l_5(x_1,x_2,x_3,x_4,x_7))+l_5(x_1,x_2,l_3(x_3,x_4,x_5),x_6,x_7)+l_5(x_1,x_2,x_5,l_3(x_3,x_4,x_6),x_7)$\\
			$+l_5(x_1,x_2,x_5,x_6,l_3(x_3,x_4,x_7))= l_3(x_3,x_4,l_5(x_1,x_2,x_5,x_6,x_7))+l_5(l_3(x_1,x_2,x_3),x_4,x_5,x_6,x_7)$\\
			$+l_5(x_3,l_3(x_1,x_2,x_4),x_5,x_6,x_7)+l_5(x_3,x_4,l_3(x_1,x_2,x_5),x_6,x_7)+l_5(x_3,x_4,x_5,l_3(x_1,x_2,x_6),x_7)$\\
			$+l_5(x_1,x_2,x_3,x_4,l_3(x_5,x_6,x_7))+l_5(x_3,x_4,x_5,x_6,l_3(x_1,x_2,x_7))$.
	\end{enumerate}
\end{definition}

A $3$-Lie $2$-algebra $(V_1,V_0,d,l_3,l_5)$ is called {\bf skeletal (strict)} if $d=0$ ($l_5=0$).

\begin{definition}\label{3lie2mor}{\rm(\cite{ZLS})}
	Let $\mathcal{V}=(V_1,V_0,d,l_3,l_5)$ and $\mathcal{V}'=(V'_1,V'_0,d',l'_3,l'_5)$ be two $3$-Lie $2$-algebras. A {\bf morphism} $\varphi$ from $\mathcal{V}$ to $\mathcal{V'}$ consists of
\begin{itemize}
\item[$\bullet$] a chain map $f:\mathcal{V}\rightarrow\mathcal{V'}$, which consists of linear maps $f_0:V_0\rightarrow V'_0$ and $f_1:V_1\rightarrow V'_1$ satisfying $f_1\circ d=d'\circ f_0$,
\item[$\bullet$] a completely skew-symmetric trilinear map $f_2:V_0\wedge V_0\wedge V_0\rightarrow V'_1$,
 \end{itemize}
such that for any $x_i\in V_0$ and $v\in V_1$, we have
	\begin{itemize}
		\item[\rm(a)] $d'(f_2(x_1,x_2,x_3))=f_0(l_3(x_1,x_2,x_3))-l'_3(f_0(x_1),f_0(x_2),f_0(x_3))$,
		\item[\rm(b)] $f_2(x_1,x_2,dv)=f_1(l_3(x_1,x_2,v))-l'_3(f_0(x_1),f_0(x_2),f_1(v))$,
		\item[\rm(c)]$l'_3(f_0(x_1),f_0(x_2),f_2(x_3,x_4,x_5))-l'_3(f_2(x_1,x_2,x_3),f_0(x_4),f_0(x_5))$\\
			$-l'_3(f_0(x_3),f_2(x_1,x_2,x_4),f_0(x_5))-l'_3(f_0(x_3),f_0(x_4),f_2(x_1,x_2,x_5))$\\
			$-l'_5(f_0(x_1),f_0(x_2),f_0(x_3),f_0(x_4),f_0(x_5))=f_2(l_3(x_1,x_2,x_3),x_4,x_5)$\\
			$+f_2(x_3,l_3(x_1,x_2,x_4),x_5)+f_2(x_3,x_4,l_3(x_1,x_2,x_5))$\\
			$-f_2(x_1,x_2,l_3(x_3,x_4,x_5))-f_1(l_5(x_1,x_2,x_3,x_4,x_5))$.
	\end{itemize}
	If $f_0$ and $f_1$ are invertible, then we call $f$ is an {\bf isomorphism}.
\end{definition}

In the following, we give the definition of $2$-derivations on $3$-Lie $2$-algebras.
\begin{definition}\label{3lie2d}
	A {\bf $2$-derivation} of a $3$-Lie $2$-algebra $\mathcal{V}$ is a triple $(X_0,X_1,l_X)$, where $X=(X_0,X_1)\in {\rm End}(V_0)\oplus {\rm End}(V_1)$ and $l_X:V_0\wedge V_0\wedge V_0\rightarrow V_1$, such that for any $x,y,z,x_i\in V_0, v\in V_1$, the following equalities hold:
	\begin{itemize}
		\item[\rm(a)]  $X_0\circ d=d\circ X_1$,
		\item[\rm(b)] $dl_X(x,y,z)\!=\!X_0l_3(x,y,z)-l_3(X_0(x),y,\!z)-l_3(x,\!X_0(y),\!z)-l_3(x,y,\!X_0(z))$,
		\item[\rm(c)] $l_X(x,\!y,\!dv)\!=\!X_1l_3(x,y,v)-l_3(X_0(x),y,\!v)-l_3(x,\!X_0(y),\!v)-l_3(x,y,\!X_1(v))$,
		\item[\rm(d)] $X_1l_5(x_1,x_2,x_3,x_4,x_5)=l_X(l_3(x_1,x_2,x_3),x_4,x_5)+l_X(x_3,l_3(x_1,x_2,x_4),x_5)$\\
		$+l_X(x_3,x_4,l_3(x_1,x_2,x_5))-l_X(x_1,x_2,l_3(x_3,x_4,x_5))+l_3(l_X(x_1,x_2,x_3),x_4,x_5)$\\
		$+l_3(x_3,l_X(x_1,x_2,x_4),x_5)+l_3(x_3,x_4,l_X(x_1,x_2,x_5))-l_3(x_1,x_2,l_X(x_3,x_4,x_5))\!$\\
		$+\sum\limits_{i=1}^5l_5(x_1,\cdots,X_0(x_i),\cdots,x_5).$
	\end{itemize}
	If $l_X=0$, then we call the triple $(X_0,X_1,0)$ is a {\bf strict $2$-derivation}  of the $3$-Lie $2$-algebra $\mathcal{V}$.
	\vspace{-0.15cm}
\end{definition}

\begin{example}
	Let $\mathcal{V}=(V_1,V_0,d,l_3,l_5)$ be a $3$-Lie $2$-algebra. Fix any $x\wedge y\in \wedge^2V_0$. Define $ad_{x\wedge y}\in{\rm End}(V_0\oplus V_1)$ by
	\begin{align*}
	ad_{x\wedge y}(z+v)=l_3(x,y,z+v),\ z\in V_0,v\in V_1.
	\end{align*}
	Set $X_0=ad_{x\wedge y}|_{V_0},X_1=ad_{x\wedge y}|_{V_1}$. Define $l_X:V_0\wedge V_0\wedge V_0\rightarrow V_1$ by
	\begin{align*}
	l_X(x_1,x_2,x_3)=-l_5(x,y,x_1,x_2,x_3), x_1,x_2,x_3\in V_0.
	\end{align*}
	Then $(X_0,X_1,l_X)$ is a $2$-derivation of $\mathcal{V}$.
\end{example}

A $3$-Lie $2$-algebra $\mathcal{V}$ with a $2$-derivation will be denoted by $(\mathcal{V};(X_0,X_1,l_X))$, which we call it a {\bf $3$-Lie$2$Der pair}. In particular, a skeletal $3$-Lie $2$-algebra (resp. a strict $3$-Lie $2$-algebra) with a $2$-derivation (resp. a strict $2$-derivation)  will be called a {\bf skeletal  $3$-Lie$2$Der pair} (resp. a {\bf strict  $3$-Lie$2$Der pair}).

\begin{definition}\label{3lie2dpiso}
	Let $(\mathcal{V};(X_0,X_1,l_X))$ and $(\mathcal{V'};(X'_0,X'_1,l'_X))$ be two {\rm 3-Lie2Der} pairs. An {\bf isomorphism} from $(\mathcal{V};(X_0,X_1,l_X))$ to $(\mathcal{V'};(X'_0,X'_1,l'_X))$ is a quadruple $(f_0,f_1,f_2,g)$, where $f_0:V_0\rightarrow V'_0$, $f_1:V_1\rightarrow V'_1$, $f_2:V_0\wedge V_0\wedge V_0\rightarrow V'_1$ and $g:V_0\rightarrow V'_1$ are linear maps, such that $(f_0,f_1,f_2)$ is a $3$-Lie $2$-algebra isomorphism and the following identities hold for all $x,y,z\in V_0$, $v\in V_1$:
	\begin{itemize}
		\item[\rm(a)]$X'_0(f_0(x))-f_0(X_0(x))=d'(g(x))$,
		\item[\rm(b)] $X'_1(f_1(v))-f_1(X_1(v))=g(d(v))$,
		\item[\rm(c)] $f_1(l_X(x,y,z))+f_2(X_0(x),y,z)+f_2(x,X_0(y),z)+f_2(x,y,X_0(z))\\
		\ \ \ -X'_1(f_2(x,y,z))-l'_X(f_0(x),f_0(y),f_0(z))=l'_3(g(x),f_0(y),f_0(z))\\
		\ \ \  +l'_3(f_0(x),g(y),f_0(z))+l'_3(f_0(x),f_0(y),g(z))-g(l_3(x,y,z)).$
	\end{itemize}
	If $f_2=g=0$, then the quadruple $(f_0,f_1,0,0)$ is called a {\bf strict isomorphism}.
\end{definition}

\begin{prop}\label{oneonesk}
	There exists a one-to-one correspondence between skeletal {\rm $3$-Lie$2$Der} pairs and triples $((L,\theta_L),(V;\rho,\theta_V),(\alpha_3,\alpha_2))$, where $(L,\theta_L)$ is a {\rm $3$-LieDer} pair, $(V;\rho,\theta_V)$ is a representation of $(L,\theta_L)$, and $(\alpha_3,\alpha_2)$ is a $3$-cocycle associated to the representation $(V;\rho,\theta_V)$.
\end{prop}
{\bf Proof.}
Suppose that $(\mathcal{V};(X_0,X_1,l_X))$ is a skeletal {\rm 3-Lie2Der} pair, where $\mathcal{V}=(V_1,V_0,0,l_3,l_5)$. By condition (d) in Definition \ref{3lie2alg}, we deduce that $V_0$ is a $3$-Lie algebra with the $3$-Lie bracket given by $l_3$. By condition (b) in Definition \ref{3lie2d}, we obtain that $X_0$ is a derivation of $V_0$ and hence $(V_0,X_0)$ is a {\rm 3-LieDer} pair.

Define $\rho:\wedge^2V_0\rightarrow {\rm End}(V_1)$ by
$$\rho(x,y)(v)=l_3(x,y,v),\quad \forall~x,y\in V_0, v\in V_1.$$
 By conditions (e) and (f) in Definition \ref{3lie2alg}, $\rho$ is a representation of $3$-Lie algebra $V_0$. Then by Definition \ref{reppa} and condition (c) in Definition \ref{3lie2d}, $(V_1;\rho,X_1)$ is a representation of {\rm $3$-LieDer} pair $(V_0,X_0)$. Moreover, by condition (g) in Definition \ref{3lie2alg}, we get ${\rm d}l_5=0$ and by condition (d) in Definition \ref{3lie2d}, we get ${\rm d}l_X-\delta l_5=0$. Thus we have
\begin{align*}
\partial(l_5,l_X)=({\rm d}l_5,{\rm d}(l_X)+(-1)^3\delta l_5)=(0,0),
\end{align*}
which implies that $(l_5,l_X)$ is a $3$-cocycle.

Conversely, let $(V;\rho,\theta_V)$ be a representation of {\rm 3-LieDer} pair $(L,\theta_L)$ and $(\alpha_3,\alpha_2)\in \mathcal{Z}_{\rm 3-LieDer}^3(L;V)$. Set $V_1=V$, $V_0=L$.
For $x,y,z\in V_0$, $u,v\in V_1$, define $l_3:V_i\times V_j\times V_k\rightarrow V_{i+j+k}$, $0\leq i+j+k\leq 1$ by
\begin{align*}
&l_3(x,y,z)=[x,y,z]_L,\\
&l_3(x,u,v)=l_3(u,x,v)=l_3(u,v,x)=0,\\
&l_3(u,x,y)=-l_3(x,u,y)=l_3(x,y,u)=\rho(x,y)(u)
\end{align*}
and for any $x_i, 1\leq i\leq 5$, define $l_5:(\wedge^2 V_0)\otimes(\wedge^3 V_0)\rightarrow V_1$ by
\begin{align*}
l_5(x_1,x_2,x_3,x_4,x_5)=\alpha_3(x_1,x_2,x_3,x_4,x_5).
\end{align*}
Then $\mathcal{V}:=(V_1,V_0,0,l_3,l_5)$ is a skeletal $3$-Lie $2$-algebra.

Define $l_X:V_0\wedge V_0\wedge V_0\rightarrow V_1$ by
\begin{align*}
l_X(x,y,z)=-l_X(y,x,z)=-l_X(x,z,y)=\alpha_2(x,y,z).
\end{align*}
Then  $(\theta_V,\theta_L,l_X)$ is a $2$-derivation of the skeletal $3$-Lie $2$-algebra $\mathcal{V}$. Therefore, $(\mathcal{V};(\theta_V,\theta_L,l_X))$ is a skeletal {\rm $3$-Lie$2$Der} pair. \hfill $\Box$
\vspace{3mm}

In the following, we shall define the equivalence relation between triples $((L,\theta_L),(V;\rho,\theta_V),(\alpha_3,\alpha_2))$ and show that there is a one-to-one correspondence between equivalent classes of such triples and isomorphism classes of skeletal {\rm $3$-Lie$2$Der} pairs.

\begin{definition}\label{3lie2dtri}
	Let $((L,\theta_L),(V;\rho,\theta_V),(\alpha_3,\alpha_2))$ and $((L',\theta_{L'}),(V';\rho',\theta_{V'}),(\alpha'_3,\alpha'_2))$ be two triples as described in Proposition \ref{oneonesk}. They are called {\bf equivalent} if there exist {\rm 3-LieDer} pair isomorphism $\sigma:(L,\theta_L)\rightarrow (L',\theta_{L'})$, linear isomorphism $\tau:V\rightarrow V'$ and two linear maps $\lambda:L\wedge L\wedge L\rightarrow V'$, $\mu:L\rightarrow V'$ such that for any $x_i,x,y,z\in L$, the following equalities hold:
	\begin{enumerate}
		\item[\rm(a)] $\theta_{V'}\circ\tau=\tau\circ\theta_V$,
		\item[\rm(b)] $\tau\circ\rho(x,y)\circ\tau^{-1}=\rho'(\sigma(x),\sigma(y))$,
		\item[\rm(c)] $\lambda([x_1,x_2,x_3]_L,x_4,x_5)+\lambda(x_3,[x_1,x_2,x_4]_L,x_5)+\lambda(x_3,x_4,[x_1,x_2,x_5]_L)\\
		-\lambda(x_1,x_2,[x_3,x_4,x_5]_L)-\tau(\alpha_3(x_1,x_2,x_3,x_4,x_5))\\
		=\rho'(\sigma(x_1),\sigma(x_2))\lambda(x_3,x_4,x_5)-\rho'(\sigma(x_3),\sigma(x_4))\lambda(x_1,x_2,x_5)\\
		-\rho'(\sigma(x_4),\sigma(x_5))\lambda(x_1,x_2,x_3)-\rho'(\sigma(x_5),\sigma(x_3))\lambda(x_1,x_2,x_4)\\
		-\alpha'_3(\sigma(x_1),\sigma(x_2),\sigma(x_3),\sigma(x_4),\sigma(x_5))$,
		\item[\rm(d)] $\tau(\alpha_2(x,y,z))+\lambda(\theta_L(x),y,z)+\lambda(x,\theta_L(y),z)+\lambda(x,y,\theta_L(z))$\\
		$-\theta_{V'}(\lambda(x,y,z))-\alpha'_2(\sigma(x),\sigma(y),\sigma(z))=\rho'(\sigma(x),\sigma(y))\mu(z)$\\
		$+\rho'(\sigma(y),\sigma(z))\mu(x)+\rho'(\sigma(z),\sigma(x))\mu(y)-\mu([x,y,z]_L).$
	\end{enumerate}

\end{definition}
\begin{theorem}\label{thm51}
	There exists a one-to-one correspondence between isomorphism classes of skeletal {\rm 3-Lie2Der} pairs and equivalent classes of triples $((L,\theta_L),(V;\rho,\theta_V),(\alpha_3,\alpha_2))$, where $(L,\theta_L)$ is a {\rm 3-LieDer} pair, $(V;\rho,\theta_V)$ is a representation of $(L,\theta_L)$, and $(\alpha_3,\alpha_2)$ is a $3$-cocycle associated to the representation $(V;\rho,\theta_V)$.
\end{theorem}
{\bf Proof.}
Let $(f_0,f_1,f_2,g)$ be an isomorphism between skeletal {\rm 3-Lie2Der} pairs $(\mathcal{V};(X_0,X_1,l_X))$ and $(\mathcal{V'};(X'_0,X'_1,l'_X))$. Let
$((V_0,X_0),(V_1;\rho,X_1),(l_5,l_X))$ and $((V'_0,X'_0),(V'_1;\rho',X'_1),(l'_5,l'_X))$ be the corresponding triples given by Proposition \ref{oneonesk}. Set $\sigma=f_0,\tau=f_1,\lambda=f_2,\mu=g$. It is obvious that $\tau:V_1\rightarrow V'_1$ is a linear isomorphism. By condition (a) in Definition \ref{3lie2mor} and condition (a) in Definition \ref{3lie2dpiso}, $\sigma:(V_0,\theta_{V_0})\rightarrow (V'_0,\theta_{V'_0})$ is a {\rm $3$-LieDer} pair isomorphism.

By conditions (b) and (c) in Definition \ref{3lie2dpiso}, we can show that conditions (a) and (d) in Definition \ref{3lie2dtri} hold.  By conditions (b) and (c) in Definition \ref{3lie2mor}, we obtain that conditions (b) and (c) in Definition \ref{3lie2dtri} hold. Thus the triples $((V_0,X_0),(V_1;\rho,X_1),(l_5,l_X))$ and $((V'_0,X'_0),(V'_1;\rho',X'_1),(l'_5,l'_X))$ are equivalent.

The converse can be proved similarly. We omit the details. \hfill $\Box$

\vspace{3mm}

In the following, we give the definition of crossed modules of {\rm $3$-LieDer} pairs.

\begin{definition}\label{cro}
	A {\bf crossed module} of {\rm 3-LieDer} pairs is a quadruple $((A,\theta_A),(B,\theta_B),(A;\rho,\theta_A),\eta)$, where $(A,\theta_A)$ and $(B,\theta_B)$ are {\rm $3$-LieDer} pairs, $(A;\rho,\theta_A)$ is a representation of {\rm $3$-LieDer} pair $(B,\theta_B)$ and $\eta:(A,\theta_A)\rightarrow(B,\theta_B)$ is a morphism of {\rm $3$-LieDer} pairs, such that for any $u,v,w\in A$, $x,y\in B$,
	\begin{align}
	\eta(\rho(x,y)(u))&=[x,y,\eta(u)]_B,\label{cro1}\\
	\rho(\eta(u),\eta(v))(w)&=[u,v,w]_A,\label{cro2}\\
	\rho(x,\eta(u))(v)&=-\rho(x,\eta(v))(u)\label{cro3},\\
\rho(x,y)[u,v,w]_A&=[\rho(x,y)(u),v,w]_A+[u,\rho(x,y)(v),w]_A+[u,v,\rho(x,y)(w)]_A.\label{cro4}
	\end{align}
\end{definition}

\begin{prop}\label{otosc}
	There exists a one-to-one correspondence between strict {\rm 3-Lie2Der} pairs and crossed modules of {\rm 3-LieDer} pairs.
\end{prop}
{\bf Proof.}
Suppose that $(\mathcal{V};(X_0,X_1,0))$ is a strict {\rm 3-Lie2Der} pair, where $\mathcal{V}=(V_1,V_0,d,l_3,0)$. Let $A=V_1$ and $B=V_0$. Define the $3$-ary bracket operations $[\cdot,\cdot,\cdot]_A$ and $[\cdot,\cdot,\cdot]_B$ by
\begin{align}
[u,v,w]_A&=l_3(du,dv,w),\label{Abrac}\\
[x,y,z]_B&=l_3(x,y,z).
\end{align}
Then it is straightforward to check that $(A,[\cdot,\cdot,\cdot]_A)$ and $(B,[\cdot,\cdot,\cdot]_B)$ are $3$-Lie algebras. Moreover, by conditions (a) and (c) in Definition \ref{3lie2d}, we have
\begin{align*}
&\ \ \ \ X_1([u,v,w]_A)=X_1l_3(du,dv,w)\\
&=l_3(X_0(du),dv,w)+l_3(du,X_0(dv),w)+l_3(du,dv,X_1(w))\\
&=l_3(d(X_1(u)),dv,w)+l_3(du,d(X_1(v)),w)+l_3(du,dv,X_1(w))\\
&=[X_1(u),v,w]_A+[u,X_1(v),w]_A+[u,v,X_1(w)]_A,
\end{align*}
which implies that $X_1$ is a derivation of $(A,[\cdot,\cdot,\cdot]_A)$. Similarly, by conditions (b) in Definition \ref{3lie2d}, we deduce that $X_0$ is a derivation of $(B,[\cdot,\cdot,\cdot]_B)$. Therefore, $(A,X_1)$ and $(B,X_0)$ are
{\rm 3-LieDer} pairs.

Define $\rho:\wedge^2B\rightarrow{\rm End}(A)$ by
\begin{align}
\rho(x,y)(u)=l_3(x,y,u).
\end{align}
By conditions (e) and (f) in Definition \ref{3lie2alg},  $\rho$ is a representation of $3$-Lie algebra $B$. By Definition \ref{reppa} and condition (c) in Definition \ref{3lie2d}, we obtain that $(A;\rho,X_1)$ is a representation of {\rm 3-LieDer} pair $(B,X_0)$.

Let $\eta=d:A\rightarrow B$. Then by condition (a) in Definition \ref{3lie2d} and the fact that
\begin{align*}
\eta([u,v,w]_A)&=dl_3(du,dv,w)=l_3(du,dv,dw)=[du,dv,dw]_B=[\eta(u),\eta(v),\eta(w)]_B.
\end{align*}
We deduce that $\eta$ is a {\rm 3-LieDer} pair morphism.

By condition (f) in Definition \ref{3lie2alg}, we have
\begin{align*}
&\ \ \ \ \rho(x,y)([u,v,w]_A)=l_3(x,y,l_3(u,dv,dw))\\
&=l_3(l_3(x,y,u),dv,dw)+l_3(u,l_3(x,y,dv),dw)+l_3(u,dv,l_3(x,y,dw))\\
&=l_3(\rho(x,y)(u),dv,dw)+l_3(u,dl_3(x,y,v),dw)+l_3(u,dv,dl_3(x,y,w))\\
&=[\rho(x,y)(u),v,w]_A+[u,\rho(x,y)(v),w]_A+[u,v,\rho(x,y)(w)]_A,
\end{align*}
which implies that \eqref{cro4} holds..

By a direct calculation, we have
\begin{align*}
\eta(\rho(x,y))(u)&=d\rho(x,y)(u)=dl_3(x,y,u)=l_3(x,y,du)=[x,y,du]_B,\\
\rho(\eta(u),\eta(v))(w)&=l_3(\eta(u),\eta(u),w)=l_3(du,dv,w)=[u,v,w]_A,\\
\rho(x,\eta(u))(v)&=l_3(x,du,v)=-l_3(x,dv,u)=-\rho(x,dv)(u),
\end{align*}
which imply that $(\ref{cro1})\!-\!(\ref{cro3})$ hold. Therefore, the quadruple $((A,X_1),(B,X_0),(A;\rho,X_1),\eta)$ is a crossed module of {\rm 3-LieDer} pairs.

The converse  can be proved similarly. We omit the details. \hfill $\Box$
\vspace{3mm}

\begin{definition}\label{croequ}
	Let $((A,\theta_A),(B,\theta_B),(A;\rho,\theta_A),\eta)$ and $((A',\theta_{A'}),(B',\theta_{B'}),(A;\rho',\theta_{A'}),\eta')$ be crossed modules of {\rm 3-LieDer} pairs. They are called {\bf equivalent} if there exist
	{\rm 3-LieDer} pair isomorphisms $\sigma:(B,\theta_B)\rightarrow(B',\theta_{B'})$ and $\tau:(A,\theta_A)\rightarrow(A',\theta_{A'})$ such that for any $x,y\in B$, the following equalities hold:
		\begin{enumerate}
		\item[\rm(a)]$\tau\circ\eta=\eta'\circ\sigma$,
		\item[\rm(b)]$\tau\circ\rho(x,y)\circ\tau^{-1}=\rho'(\sigma(x),\sigma(y)).$
	\end{enumerate}
\end{definition}

\begin{theorem}\label{thm52}
	There exists a one-to-one correspondence between strict isomorphism classes of strict {\rm 3-Lie2Der} pairs and equivalent classes of crossed modules of {\rm 3-LieDer} pairs.
\end{theorem}
{\bf Proof.}
Let $(f_0,f_1,0,0)$ be the  strict isomorphism  between the strict {\rm 3-Lie2Der} pairs $(\mathcal{V};(X_0,X_1,0))$ and $(\mathcal{V'};(X'_0,X'_1,0))$.   Let $((V_1,X_1),(V_0,X_0),(\rho,\!V_1,\!X_1),\eta)$ and $((V'_1,X'_1),(V'_0,X'_0),(\rho',\!V'_1,\!X'_1),\eta')$ be crossed modules of {\rm 3-LieDer} pairs, corresponding to $(\mathcal{V};(X_0,X_1,0))$ and $(\mathcal{V'};(X'_0,X'_1,0))$ given by Proposition \ref{otosc} respectively. By conditions (a) and (b) in Definition \ref{3lie2dpiso} , we get that $\sigma:(V_0,X_0)\rightarrow(V'_0,X'_0)$ and $\tau:(V_1,X_1)\rightarrow(V'_1,X'_1)$ are {\rm 3-LieDer} pair isomorphisms.

Set $\sigma=f_0$, $\tau=f_1$, $\eta=d$ and $\eta'=d'$.  Define $\rho:\wedge^2V_0\rightarrow {\rm End}(V_1)$ and  $\rho':\wedge^2V'_0\rightarrow {\rm End}(V'_1)$ by
\begin{eqnarray*}
\rho(x,y)(u)&=&l_3(x,y,u),\ x,y\in V_0,\ u\in V_1,\\
\rho(x',y')(u')&=&l'_3(x',y',u'),\ x',y'\in V'_0,\ u'\in V'_1	.
\end{eqnarray*}
By the condition $f_1\circ d=d'\circ f_0$ in Definition \ref{3lie2mor}, we have $\tau\circ\eta=\eta'\circ\sigma$. By condition (b) in Definition \ref{3lie2mor}, we deduce that $\tau\circ\rho(x,y)\circ\tau^{-1}=\rho'(\sigma(x),\sigma(y)).$ Therefore, the above two crossed modules of {\rm 3-LieDer} pairs are equivalent.

The converse  can be proved similarly. We omit the details. \hfill $\Box$

\vspace{3mm}
\noindent
{\bf Acknowledgements. } This research is supported by NSFC (11901501).

 \end{document}